\documentclass{amsart} 
\usepackage{amssymb}
\usepackage{amsmath}
\usepackage{amsfonts}

\begin{document}
\newtheorem{theo}{Theorem}[section]
\newtheorem{prop}[theo]{Proposition}
\newtheorem{lemma}[theo]{Lemma}
\newtheorem{exam}[theo]{Example}
\newtheorem{coro}[theo]{Corollary}
\theoremstyle{definition}
\newtheorem{defi}[theo]{Definition}
\newtheorem{rem}[theo]{Remark}


\newcommand{\Ac}{{\mathcal A}}
\newcommand{\Bc}{{\mathcal B}}
\newcommand{\Cc}{{\mathcal C}}
\newcommand{\Dc}{{\mathcal D}}
\newcommand{\Ic}{{\mathcal I}}
\newcommand{\Jc}{{\mathcal J}}
\newcommand{\Lc}{{\mathcal L}}
\newcommand{\Oc}{{\mathcal O}}
\newcommand{\Pc}{{\mathcal P}}
\newcommand{\Sc}{{\mathcal S}}
\newcommand{\Uc}{{\mathcal U}}
\newcommand{\Zi}{{\rm Z}^i_1}

\newcommand{\ax}{{\rm ax}}
\newcommand{\Acc}{{\rm Acc}}
\newcommand{\Act}{{\rm Act}}
\newcommand{\Bd}{{\rm Bd}}
\newcommand{\ded}{{\rm ded}}
\newcommand{\Gm}{{$\Gamma_0$}}
\newcommand{\ID}{{${\rm ID}_1^i(\Oc)$}}
\newcommand{\PAP}{{${\rm PA}(P)$}}
\newcommand{\ACA}{{${\rm ACA}^i$}}
\newcommand{\Rd}{{\rm Rd}}
\newcommand{\RefP}{{${\rm Ref}^*({\rm PA}(P))$}}
\newcommand{\RefS}{{${\rm Ref}^*({\rm S}(P))$}}
\newcommand{\Rfn}{{\rm Rfn}}
\newcommand{\Seq}{{\rm Seq}}
\newcommand{\tar}{{\rm Tarski}}
\newcommand{\Typ}{{\rm Typ}}
\newcommand{\UNFA}{{${\mathcal U}({\rm NFA})$}}

\author{Nik Weaver}

\title [Predicativity beyond \Gm{}]
       {Predicativity beyond \Gm{}}

\address {Department of Mathematics\\
          Washington University in Saint Louis\\
          Saint Louis, MO 63130}

\email {nweaver@math.wustl.edu}

\date{\em January 13, 2007}

\begin{abstract}
We reevaluate the claim that predicative reasoning (given the natural
numbers) is limited by the Feferman-Sch\"utte ordinal \Gm{}. First we
comprehensively criticize the arguments that have been offered in support
of this position. Then we analyze predicativism from first principles and
develop a general method for accessing ordinals which is predicatively
valid according to this analysis. We find that the Veblen ordinal
$\phi_{\Omega^\omega}(0)$, and larger ordinals, are predicatively provable.
\end{abstract}

\maketitle


The precise delineation of the extent of predicative reasoning is
possibly one of the most remarkable modern results in the foundations
of mathematics. Building on ideas of Kreisel \cite{Kre1, Kre2},
Feferman \cite{Fef1} and Sch\"utte \cite{Sch1, Sch2} independently
identified a countable ordinal \Gm{} and argued that it is the
smallest predicatively non-provable ordinal. (Throughout, I take
``predicative'' to mean ``predicative given the natural numbers''.)
This conclusion has become the received view in the foundations community,
with reference \cite{Fef1} in particular having been cited with approval
in virtually every discussion of predicativism for the past forty years.
\Gm{} is now commonly referred to as ``the ordinal of predicativity''.
Some recent publications which explicitly make this assertion are
\cite{Ant, Ara, Avi0, Avi1, AviS, Fef9, FefSt, Fri, Hel, JKS, Kah,
Poh, Rat, Rue, Sim2, Str}.

This achievement is notable both for its technical sophistication and
for the insight it provides into an important foundational stance.
Although predicativism is out of favor now, at one time it was
advocated by such luminaries as Poincar\'e, Russell, and Weyl.
(Historical overviews are given in \cite{Fef9} and \cite{Par}.)
Its central principle --- that sets have to be ``built up from
below'' --- is, on its face, reasonable and attractive. With its
rejection of a metaphysical set concept, predicativism also
provides a cogent resolution of the set-theoretic paradoxes and
is more in line with the positivistic aspect of modern analytic
philosophy than are the essentially platonic views which have become
mathematically dominant.

Undoubtedly one of the main reasons predicativism was not accepted
by the general mathematical public early on was its apparent failure
to support large portions of mainstream mathematics. However, we now
know that the bulk of core mathematics can in fact be developed in
predicative systems \cite{Fef5, Sim}, and the limitation identified
by Feferman and Sch\"utte is probably now a primary reason, possibly
{\it the} primary reason, for predicativism's nearly universal
unpopularity.${}^1$ There do exist important mainstream theorems
which are known to in various senses require
provability of \Gm{}, and in any case \Gm{} is sufficiently tame that it
is simply hard to take seriously any approach to foundations that prevents
one from recognizing ordinals at least this large. Thus, it is of great
foundational interest to examine carefully whether the \Gm{} limitation
really is correct. If it is not, predicativism could be more viable than
previously thought and its current peripheral status in the philosophy of
mathematics may need to be reconsidered.

I believe \Gm{} has nothing to do with predicativism. I will argue that
the current understanding of predicativism is fundamentally flawed and
that a more careful analysis shows the ``small'' Veblen ordinal
$\phi_{\Omega^\omega}(0)$, and probably much larger ordinals, to be
within the scope of predicative mathematics. It is my hope that this
conclusion will open the way to a serious reappraisal of the significance
and interest of predicativism. Elsewhere I introduce the term ``mathematical
conceptualism'' for the brand of predicativism considered here and make a
case that it is cogent, rigorous, attractive, and better suited to ordinary
mathematical practice than all other foundational stances \cite{W1}.

\section {A critique of the \Gm{} thesis}\label{sect1}

At issue is the assertion that there are well-ordered sets of all
order types less than \Gm{} and of no order types greater
than or equal to \Gm{} which can be proven to be well-ordered
using predicative methods (cf.\ \cite{Fef1}, p.\ 13 or \cite{Sch3}, p.\
220). I call this {\it the \Gm{} thesis}.

As stated, this claim is imprecise because the classical concept of
well-ordering has a variety of formulations which are not predicatively
equivalent (see \S \ref{sect1d} and \S \ref{sect2d}). In fact, previous
discussions of predicativism have tended to ignore this distinction, and
this will emerge as a crucial source of confusion (see \S \ref{sect1d}).
To fix ideas I will use the term ``well-ordered''
to mean of a set $X$ that it is totally ordered and if $Y \subseteq X$
is progressive then $Y = X$. {\it Progressive} means that for every
$a \in X$, if $\{b \in X: b \prec a\} \subseteq Y$ then $a \in Y$.

In principle, to falsify the \Gm{} thesis I need only produce (1) a
well-ordering proof of an ordered set that is isomorphic to
\Gm{} and (2) a convincing case that the proof is predicatively
valid. However, no matter how convincing I could make that case,
in light of the broad and sustained acceptance the thesis
has enjoyed it would be unsatisfying to leave the matter there.
The \Gm{} thesis has been repeatedly and forcefully defended by
two major figures, Feferman and Kreisel. Many current authors
simply assert it as a known fact. The only substantial published
criticism of which I am aware appears in \cite{How}, but even that
is somewhat ambivalent and seems to conclude in favor of the thesis.
Therefore, I take it that I have a burden not only to positively
demonstrate the power of predicative reasoning, but also to show
where the more pessimistic previous assessments went wrong.

This is a somewhat lengthy task because a great deal has been written
in support of the \Gm{} thesis from a variety of points of view.
However, I believe that the entire body of argument is specious and
can be decisively refuted. The goal of Section \ref{sect1} is to do
this in some detail.

One point before I begin. Predicativism is a philosophical position,
and prior to the acceptance of a particular formalization there will
be room for argument over its precise nature. Thus, a debate about
the \Gm{} thesis could easily degenerate into a
purely semantic dispute over the meaning of the term ``predicative''.
I therefore want to emphasize that the central claim of this section
is that there is {\it no} coherent philosophical stance which would
lead one to accept every ordinal less than \Gm{} but not \Gm{}
itself. My polemical technique will be to examine various formal
systems that have been alleged to model predicative reasoning and
indicate how in each instance the informal principles that motivate
the system actually justify a stronger system which goes beyond \Gm{}.
This is obviously independent of any special views one may have about
predicativity. (A second major claim is that none of these
supposedly predicative systems is actually predicatively legitimate.
Evaluating the justice of this claim does require some understanding
of predicativism, and I refer the reader who is not prepared to accept
any assertion of this type at face value to Section \ref{sect2}, where
I develop my views on predicativism in detail.)

\subsection{Formal systems for predicativity}\label{sect1a}

A variety of formal systems have been proposed as modelling predicative
reasoning in some form of second order arithmetic. Among the
main examples are $\Sigma$ \cite{Kre1}, $H^+$, $R^+$, $H$, $R$
\cite{Fef1}, $RA^*$ \cite{Sch3}, $P + \exists/P$ \cite{Fef3},
\RefP{} \cite{Fef4}, and \UNFA{} \cite{FefSt}. I give here a
brief sketch of their most important features.

The systems $\Sigma$, $H$, and $R$ are similar in broad outline and
need not be distinguished in this discussion; likewise for the systems
$H^+$, $R^+$, and $RA^*$.  All six express a concept of ``autonomy''
which allows one to access larger ordinal notations from smaller ones
and I will refer to them generally as ``autonomous systems''.
In the last three systems the idea is that once a predicativist
has proven the well-foundedness of a set of order type $\alpha$
he is allowed to use infinite proof trees of height $\alpha$ to
establish the well-foundedness of larger order types,
the key infinitary feature being an ``$\omega$-rule''
which permits deduction of the formula $(\forall n)\, \Ac(n)$,
where $n$ is a number variable, from the family of formulas
$\Ac(\overline{n})$ with $\overline{n}$ ranging over all numerals.
In the first three systems all proofs are finite and the key proof
principle is a ``formalized $\omega$-rule'' schema which, for each
formula $\Ac$, concludes the formula $(\forall n)\, \Ac(n)$
from a premise which arithmetically expresses that for every
number $n$ there is a proof of $\Ac(\overline{n})$. This leads to a
hierarchy of systems $S_a$ where $a$ is an ordinal notation and
$S_{a\oplus 1}$ incorporates a formalized $\omega$-rule schema
referring to proofs in $S_a$. The predicativist is then permitted
to execute a finite succession of proofs in various $S_a$'s,
subject to the requirement that passage to any $S_a$ must be
preceded by a proof that $a$ is an ordinal notation.

The linked systems $P$ and $\exists/P$ are notable for their
conception of predicativists as having a highly restricted yet not
completely trivial ability to deal with second order quantification, in
particular their being able to use only free or, to a limited extent,
existentially quantified set variables. Second order existential
quantification is actually permitted only in the ``auxiliary'' system
$\exists/P$, but once a functional has been shown to exist uniquely
one is allowed to introduce a symbol for it which can then be used
in $P$. By passing back and forth between $P$ and $\exists/P$ one is
able to produce functionals which provably enumerate larger and larger
initial segments of the ramified hierarchy over an arbitrary set
and use them to prove the well-ordering property for successively
larger ordinal notations.

The system \RefP{} is obtained by applying a general construction
${\rm Ref}^*$ to a ``schematic'' form \PAP{} of Peano arithmetic.
This construction involves extending the language of \PAP{} to allow
assertions of truth and falsehood and adding axioms
which govern the use of the truth and falsehood predicates. Paradoxes
arising from a self-referential notion of truth are avoided by regarding
the truth and falsehood predicates as partial and axiomatizing them in a
way that expresses their ultimate groundedness in facts about \PAP{}.
The ability to reason about truth in effect implements the formalized
$\omega$-rule mentioned above, and this again enables one to prove
the well-foundedness of successively larger ordinal notations. The
general idea is that \RefS{} embodies what one ``ought to accept''
given that one accepts a schematic theory ${\rm S}(P)$, and an
argument can then be made that predicativism is fundamentally based
on Peano arithmetic and therefore \RefP{} precisely captures what a
predicativist ought to accept.

Like \RefP{}, \UNFA{} is an instance of a general construction which
applies to any schematic formal system and is supposed to embody what
one ought to accept once one accepts that system. However, the exact
claim is slightly different: here we are concerned with
determining ``which operations and predicates, and which principles
concerning them, ought to be accepted'' once one has accepted the
initial system (\cite{FefSt}, p.\ 75). This problem is approached
from the point of view of generalized recursion theory and one is
allowed to generate operations and predicates by using a least fixed
point operator. It is easily seen that this recursive generation
procedure rapidly recovers Peano arithmetic from a weaker theory
${\rm NFA}$ (``non-finitist arithmetic''); therefore, \UNFA{} is already
supposed to capture predicative reasoning.

\subsection{Outline of the critique.}\label{sect1b}

All of the proposed formalizations of predicative reasoning cited
in \S \ref{sect1a} have the same provable ordinals, namely all
ordinals less than \Gm{}. This by itself might be seen as good
evidence in favor of the \Gm{} thesis simply because it seems
unlikely that so many different approaches should all have settled
on the same wrong answer.

Nonetheless, each of the formal systems in \S \ref{sect1a} is
simultaneously too weak and too strong to faithfully model
predicative reasoning and thereby verify the claim about
\Gm{}. They are all too weak for a general reason I discuss in
\S \ref{sect1c}; in brief, anyone who accepts a given system
ought to be able to grasp its global validity and then go beyond it.
This is an old objection and there are several responses to it on
record. However, these responses, which I review below, are not
well-taken because they typically involve postulating (usually with
little or no overt justification) a limitation on predicative reasoning
which, if true, would actually have prevented a predicativist
from working within the original system.

In addition, each system is manifestly impredicative in some way, and
hence too strong. This fact does not seem to be widely appreciated,
but it is hardly obscure. The autonomous systems impredicatively
infer a transfinite iteration of reflection principles from a statement
of transfinite induction. $P + \exists/P$ allows predicate substitution
for $\Sigma^1_1$ formulas, so that for every $\Sigma^1_1$ formula $\Ac$
it in effect lets one reason about $\{n: \Ac(n)\}$ as if this were a
meaningful set, which in general is predicatively not the case. \RefP{}
makes truth claims about schematic predicates which do not make sense
unless one assumes an impredicative comprehension axiom. \UNFA{}
employs a patently impredicative least fixed point operator and also
treats schematic predicates in a way that again can only be justified
by impredicative comprehension. I will elaborate on all of these points
below.

The most striking impredicativity is the least fixed point operator of
\UNFA{}, but the other instances are actually more significant because
they fit into a general pattern. The basic problem is that in each of
these systems one proves the well-foundedness of successively larger
ordinal notations by an inductive argument that at each step involves
generating some kind of iterative hierarchy which is used to prove transfinite
induction at the next level --- but this does not justify the statement
of transfinite {\it recursion} which is needed to generate the next
hierarchy. In order to make this inference from induction to recursion
one has to smuggle an impredicative step somewhere into the proof, and
this is the function of all the other examples of impredicativity noted
above. I will return to this point in \S \ref{sect1i}.

There are also more subtle instances of impredicativity which occur in
the use of self-applicative schematic predicates in \RefP{} and \UNFA{}
and in the use of self-applicative truth and falsehood predicates in
\RefP{}; see \S \ref{sect2e} and \S \ref{sect2f}.

I next describe an objection that is generally applicable, and then
go on to discuss the individual formal systems.

\subsection{A general difficulty}\label{sect1c}

Suppose $A$ is a rational actor who has adopted some foundational
stance. Any attempt to precisely characterize the limits of $A$'s
reasoning must meet the following objection: if we could show that
$A$ would accept every member of some set of statements $\Sc$, then
$A$ should see this too and then be able to go beyond $\Sc$, e.g.\
by asserting its consistency. Thus, $\Sc$ could not have been a
complete collection of all the statements (in a given language)
that $A$ would accept. A similar argument can be made about
attempts to characterize $A$'s provable ordinals.

There are a variety of ways in which this objection might be overcome.
$A$ may actually be unable to recognize $\Sc$ as a legitimate set, for
instance if $\Sc$ is infinite and $A$ is a finitist. Or the language
in use may not be capable of expressing the consistency of $\Sc$. Or
perhaps $A$ can indeed see, as we do, that there exists a proof that
he would accept for each statement in $\Sc$, but he cannot go from
this to actually accepting every statement in $\Sc$ (though it is
difficult to imagine a plausible set of beliefs that would not allow
him to take this step). Or it may be possible to identify some special
limitation in $A$'s belief system which prevents him from grasping the
validity of all of $\Sc$ at once despite his ability to accept each
statement in $\Sc$ individually.${}^2$

Defenses of the \Gm{} thesis generally take the last approach. This is
tricky for a slightly subtle reason. It is not hard to believe that
$A$ (or anyone) is unable to simultaneously identify exactly which
statements are true from his perspective. But it is more difficult
to reconcile this with the claim that we {\it do} know that he would
accept each statement in $\Sc$. The most obvious way to establish this
claim would be to explicitly show how $A$ would prove each statement in
$\Sc$, and this is actually the method used in the case at hand: each
of the proposed formal systems for predicativism is accompanied by a
recursive proof schema which is supposed to show how a predicativist
could use the system to access every ordinal less than \Gm{}. What
is confusing here is the suggestion that we can see that he would
accept each proof in the schema but he cannot see this.

In fact this is highly implausible, for the following reason.
Let $\gamma_0 = 1$ and $\gamma_{n+1} = \phi_{\gamma_n}(0)$, where
$\{\phi_\alpha(\beta)\}$ is the Veblen hierarchy of critical
functions, so that $\Gamma_0 = \sup_{n \in \omega} \gamma_n$. Now
in general we are not merely given a recursive set of proofs which
establish for each $n$ that some notation $a_n$ for $\gamma_n$ is
an ordinal notation; for each of the formal systems under discussion,
at least at an intuitive level we have a {\it single} proof that
for any $n$, $a_n$ is an ordinal notation implies $a_{n+1}$ is an
ordinal notation. It therefore becomes hard to believe that someone
who is presumed to grasp induction on $\omega$ (and even, allegedly,
in ``schematic'' form \cite{Fef4, Fef6, FefSt}) would not be able to
infer the single assertion that $a_n$ is an ordinal notation for all $n$.

It is reasonable to expect that if a predicativist understands how to go
from $a_n$ to $a_{n+1}$ for any single value of $n$, and if the passage
from $a_n$ to $a_{n+1}$ is essentially the same for all $n$, then he
{\it can} infer the statement that every $a_n$ is an ordinal notation.
As this would enable him to immediately deduce the well-foundedness of
an ordered set isomorphic to $\Gamma_0$, advocates of the $\Gamma_0$
thesis have a crucial burden to explain why he cannot in fact do this.
Yet the handful of attempts to establish this point that appear in the
literature are brief, vague, and, I will argue, simply unpersuasive.

I now turn to the systems introduced in \S \ref{sect1a}.

\subsection{The finitary autonomous systems}\label{sect1d}

The initial idea behind the finitary autonomous systems $\Sigma$,
$H$, and $R$ is that if a predicativist trusts some formal system
for second order arithmetic, say ${\rm ACA}_0$ (see \cite{Sim}),
then he should accept not only the theorems of the system itself,
but also additional statements such as ${\rm Con}({\rm ACA}_0)$
which reflect the fact that the axioms are true. Feferman \cite{Fef0}
analyzed several such ``reflection principles'' and found the strongest
of them to be the formalized $\omega$-rule schema
$$(\forall n)\left[{\rm Prov}(\ulcorner{\Ac(\overline{n})}\urcorner)
\rightarrow \Ac(n)\right],$$
where $\ulcorner{\Ac(\overline{n})}\urcorner$ is the G\"odel number of
$\Ac(\overline{n})$ and ${\rm Prov}$ formalizes ``is the G\"odel number
of a provable formula'' (here, provable in ${\rm ACA}_0$).

Having accepted this schema, the argument runs, the predicativist
is then committed to a stronger system consisting of ${\rm ACA}_0$
plus the $\omega$-rule schema, and he should therefore now accept a
version of the formalized $\omega$-rule schema which refers to provability
in this stronger system. This process can be transfinitely iterated,
yielding a family of formal systems $S_a$ indexed by Church-Kleene ordinal
notations $a$. Kreisel's idea \cite{Kre1} was that a predicativist
should accept the system indexed by $a$ when and only when he has a
prior proof that $a$ is an ordinal notation.${}^3$

Feferman \cite{Fef1} proved that when this procedure is carried out
starting with a reasonable base system $S_0$, \Gm{} is the smallest
ordinal with the property that there is no finite sequence
of ordinal notations $a_1, \ldots, a_n$ with $a_1$ a notation for
$0$, $a_n$ a notation for \Gm{}, and such that $S_{a_i}$ proves
that $a_{i+1}$ is an ordinal notation ($1 \leq i < n$). Thus, \Gm{}
is the smallest predicatively non-provable ordinal.

There are two fundamental problems with this analysis. The first is
that the plausibility of inferring soundness of $S_a$ from the fact
that $a$ is an ordinal notation hinges on our conflating two versions
of the concept ``ordinal notation'' --- supports transfinite induction
for arbitrary sets versus supports transfinite induction for arbitrary
properties --- which are {\it not} predicatively equivalent. What we
actually prove about $a$ is that, for a given partial order $\prec$
on a subset of $\omega$, if $X$ is a set with the property that
$$(\forall b)\left[(\forall c \prec b)(c \in X)
\rightarrow b \in X\right]$$
then every $b \prec a$ must belong to $X$. Classically this entails
that for every formula $\Ac$ the statement
$$(\forall b)\left[(\forall c \prec b)\, \Ac(c)
\rightarrow \Ac(b)\right]$$
implies $\Ac(b)$ for all $b \prec a$ because we can use a comprehension
axiom and reason about the set $X = \{b: \Ac(b)\}$. Predicatively
this should still be possible if, for example, $\Ac$ is arithmetical,
but not in general. Now the statement $\Pc(b) \equiv$
``if ${\rm Prov}_{S_b}(\ulcorner{\Ac}\urcorner)$ then $\Ac$, for
every formula $\Ac$'' is not only not arithmetical, it cannot even
be formalized in the language of second order arithmetic. So we
should not expect there to be any obvious way to {\it predicatively}
infer $\Pc(a)$ from what we have proven about $a$. Indeed, there are
good reasons to suppose that this inference is not legitimate, for
instance the fact that $S_a$ proves the existence of arithmetical
jump hierarchies up to $a$, which is formally stronger than the fact
that transfinite induction holds up to $a$ for sets.

One may be tempted to dismiss this first objection as technical and to grant
that predicativists can make the disputed inference, but that leads to
a second basic problem: if a predicativist could somehow infer
the soundness of $S_a$ then he actually ought to be able to infer more.
This point was made well by Howard \cite{How}. I would put it this way:
according to Kreisel, a predicativist is (somehow) always able to make
the deduction
$${\rm from }\quad I(\overline{a})\quad{\rm and}\quad
{\rm Prov}_{S_a}(\ulcorner{\Ac(\overline{n})}\urcorner),
\qquad{\rm infer}\quad \Ac(\overline{n}),\eqno{(*)}$$
where $I(a)$ formalizes the assertion that $a$ is an ordinal notation.
Shouldn't he then accept the assertion
$$(\forall a)(\forall n)\left[ I(a) \wedge
{\rm Prov}_{S_a}(\ulcorner{\Ac(\overline{n})}\urcorner)
\rightarrow \Ac(n)\right]\eqno{(**)}$$
for any formula $\Ac$?

As a straightforward consequence of \cite{Fef1}, one can use ($**$)
to prove $I(\overline{a})$ with $a$ some standard notation for \Gm{}.${}^4$
The claim must therefore be that a predicativist can recognize each
instance of ($*$) to be valid but cannot recognize the validity of
the general assertion ($**$). In other words, whenever he has proven
that $a$ is an ordinal notation he can infer the statement that all
theorems of $S_a$ hold, but he does {\it not} accept the statement
``if $a$ is an ordinal notation then all theorems of $S_a$ hold.''
Why not?
\medskip

\noindent {\bf (a)} {\it Kreisel's first answer.}
Kreisel addresses this point in \cite{Kre1}. He writes:

\begin{quotation}
\noindent Here, too, though each extension is predicative provided $<$
has been recognized by predicative means to be a well-ordering, the
general extension principle is not since [it is framed in terms of]
the concept of predicative proof [which] has no place in predicative
mathematics. (\cite{Kre1}, p.\ 297; see also p.\ 290)
\end{quotation}

Although this comment sounds authoritative, it does not hold up
under scrutiny because in whatever sense it could be said that
($**$) presumes the concept of predicative proof, the same is
true of any instance of ($*$). If we had no concept of proof or
validity then we ought not be able to make the inference ($*$)
in any instance. One can try to read something more subtle into
Kreisel's comment, but I have not found any way to elaborate it
into a convincing argument. Perhaps the best attempt appears in
\S \ref{sect1d} (b) below.

Similar reasoning would actually support a more severe conclusion.
Consider:

\begin{quotation}
\noindent Although the inference of $\Bc$ from $\Ac$ is predicative provided
$\Ac$ has been recognized by predicative means to imply $\Bc$, the general
principle of modus ponens is not since it is framed in terms of the
concept of predicative truth, which has no place in predicative
mathematics.
\end{quotation}

\noindent This is a parody, but not a gross one. In fact, I do not
really see what could make one accept the first statement and not the
second. (The rejoinder that modus ponens is not framed
specifically in terms of {\it predicative} truth misses the
point. To a predicativist, ``truth'' and ``predicative
truth'' are the same thing, so it would not make sense to suggest
that he can reason about truth but not about predicative truth.)
If one did accept the second statement, of course, this would
prevent any use of reflection principles since absent a general
grasp of modus ponens the mere acceptance of a set of axioms
would not entitle one to globally infer the truth of all
theorems provable from those axioms.
\medskip

\noindent {\bf (b)} {\it Kreisel's second answer.}
A second argument in response to something like the objection raised
above was made by Kreisel (\cite{Kre5}, \S 3.631) and cited with
approval by Feferman (\cite{Fef2}, p.\ 134). Unfortunately, the cited
passage is rather inscrutable, so it is hard to be sure what
Kreisel had in mind. I think it is something like this.
Predicativists are at any given moment only able to reason about
those subsets of $\omega$ that have previously been shown to exist.
The ``basic step'' of predicative reasoning is thus the passage from
one level $N_\alpha$ of the ramified hierarchy over $\omega$ to the
next ($N_{\alpha + 1} =$ the subsets of $\omega$ definable by second
order formulas relativized to $N_\alpha$). Now the proof that (a notation
for) $\gamma_{n+1}$ is well-founded uses only sets in $N_{\gamma_n}$, so
once $N_{\gamma_n}$ is available this proof can be executed and one can
pass to $N_{\gamma_{n+1}}$. However, we cannot go directly from
$N_{\gamma_n}$ to $N_{\gamma_{n+2}}$ since the proof that $\gamma_{n+2}$
is well-founded uses sets in $N_{\gamma_{n+1}}$ which are not yet
available. Thus, we cannot grasp the validity of the sequence of
proofs as a whole since later proofs involve the use of sets that
are not known to exist at earlier stages. Each individual proof is
admissible, however, since there is a finite stage in the reasoning
process at which the sets needed for that proof become available.

This neatly answers the question raised in \S \ref{sect1c} as to
how each proof could be recognized as valid while the entire sequence of
proofs cannot. But wait. Exactly how would one use the well-foundedness
of $\gamma_{n+1}$ proven using sets in $N_{\gamma_n}$ to ``pass to
$N_{\gamma_{n+1}}$'' and make sets at that stage available for future
proofs? If we accept Kreisel's premise then it would seem that we cannot
directly go even two levels up from $N_{\gamma_n}$ to $N_{\gamma_n + 2}$,
let alone all the way to $N_{\gamma_{n+1}}$, because the construction of
$N_{\gamma_n + 2}$ uses sets in $N_{\gamma_n+1}$ which are not yet available.
Thus, the argument that prevents us from getting up to
$\Gamma_0$ should be equally effective at preventing us from getting
from $\gamma_n$ to $\gamma_{n+1}$.

This point may become clearer if we ask how a predicativist could
establish the existence of $N_\omega$. Starting with $N_0 = \emptyset$,
he can use the basic step to
directly pass to $N_1$, then to $N_2$, and so on, so that for each
$n \in \omega$ he can give a finite proof of the existence of $N_n$. But
in order to accept the existence of $N_\omega$ he has to somehow globally
grasp that $N_n$ exists for all $n$ {\it without} sequentially proving
their existence one at a time. Presumably he can accomplish this by
recognizing the general principle that the existence of $N_{n+1}$ follows
from the existence of $N_n$ and then making an induction argument. So
evidently in this case he {\it can} accept the validity of the sequence
of proofs as a whole despite the fact that later proofs involve the use
of sets that are not known to exist at earlier stages. That is, {\it just
getting up to $N_\omega$ already requires some ability to reason
hypothetically about sets that are not yet available}. So Kreisel's
argument (if this really is what he meant) appears to make little sense.

However, this entire discussion is speculative until we are told
precisely why the proof that $\gamma_{n+1}$ is well-founded is supposed
to legitimate passage to $N_{\gamma_{n+1}}$. This takes us to Kreisel's
final argument.
\medskip

\noindent {\bf (c)} {\it Kreisel's third answer.}
Kreisel's most sophisticated analysis appears in \cite{Kre6}. Here he
rightly addresses the central question of exactly how a predicativist
would infer soundness of $S_a$ once $I(a)$ has been proven. On my reading,
the novel idea is that this inference (or something like it) would not be
based on genuinely ``understanding'' the well-ordering property of $a$,
which he now denies a predicativist could do, but instead would be directly
extracted from the structure of the proof of $I(a)$. If $P$ is the property
``the (formal) definitions at a level of the [ramified] hierarchy considered
are understood if our basic concepts are understood'' (\cite{Kre6}, p.\ 498),
then

\begin{quotation}
\noindent Since we do not have an explicit definition for $P$ $\ldots$
it seems reasonable to suppose that the {\it formal derivation of the
well-foundedness of $\beta$ is needed} $\ldots$ specifically, we expect
{\it to use the derivation as a} (naturally, infinite) {\it schema}
which need be applied only to instances of $P$ whose meaning {\it is}
determined at stage $\alpha$. (\cite{Kre6}, pp.\ 498-499; italics in original)
\end{quotation}

\noindent He adds in a footnote: ``It seems likely that the work of
Feferman and Sch\"utte `contains' all the formal details needed
$\ldots$ the principal problem is conceptual: to formulate properly
just what details are needed.''

It seems even more reasonable to suppose that if, ten years after his
first attempt (in \cite{Kre1}) to refute the objection about ($**$),
Kreisel is still not sure how to do this, then the objection is probably
valid. Here he gives us not a fully realized refutation, but merely a
speculation as to how one might be obtained. I do not think any attempt
of this type is likely to succeed for the reasons discussed at the
beginning of this section, in particular the fact that $S_a$ proves
the existence of arithmetical jump hierarchies up to $a$ and this
seems not to be predicatively entailed by $I(a)$ (cf.\ the end of $\S$
\ref{sect2d}). Moreover, even if
one could work out some way of converting formal derivations of
well-ordering in the autonomous systems into informal verifications
of soundness in some metatheory, then presumably the metatheory and
the conversion process could be formalized, and then a predicativist
should be able to apply a single instance of the formalized
$\omega$-rule to the metatheory and deduce ($**$) as a general
principle. But again, this discussion is hypothetical.
\medskip

\noindent {\bf (d)} {\it Feferman's position.}
In \cite{Fef3} Feferman raises a version of the objection and notes that
it ``involve[s] the ordinal character of the proposal via progressions,
hence [does] not apply to $P$ [$+ \exists/P$]'' (\cite{Fef3}, p.\ 85).
Similar comments appear in (\cite{Fef4}, p.\ 3) and (\cite{Fef9}, p.\ 24).
It is certainly true that the systems $P + \exists/P$, \RefP{}, and \UNFA{}
do not presume any special ability to reason using well-ordered
sets. However, Feferman nowhere openly repudiates the earlier systems, and
I read his remark in \cite{Fef4} as implying that the later systems are merely
more ``perspicuous'' than the earlier ones because they do not assume that
predicativists have any understanding of ordinals. As far as I know he has
never addressed the argument that a grasp of ordinals sufficient to justify
($*$) would also justify ($**$) and hence lead one beyond \Gm{}.${}^5$

\subsection{The infinitary autonomous systems}\label{sect1e}

In order to evaluate the infinitary (semiformal) autonomous systems we
must first clarify in exactly what way these systems are supposed to
model predicative reasoning. Surely they are not meant to be taken
literally in this regard. Perhaps we can conceive of an idealized
predicativist living in an imaginary world who is capable of actually
executing proofs of transfinite length, but in this case the allowed
proof lengths would merely depend on the nature of the imagined world,
not on which well-ordering statements the predicativist is able to prove.

Presumably the infinitary autonomous systems are meant to be taken
as modelling what an actual predicativist would consider a valid but
idealized reasoning process. In other words, the predicativist does
not himself reason within any of these systems, but he believes that
in principle they would prove true theorems if they could somehow be
implemented (in some imaginary world). On this interpretation the fact
that an ordinal $\alpha$ is autonomous within one of these systems could
lead a predicativist to accept the well-foundedness of (some notation for)
$\alpha$ {\it only if he knew this fact}. However, the only way he could
know which ordinals are autonomous is via some kind of meta-argument about
what is provable in the given system. This immediately suggests that
he should be able to get beyond \Gm{} by performing a single act of
reflection on the finitary (formal) system in which he actually reasons.

We can now see that just as in the case of the finitary autonomous systems
one is faced with a dilemma: (1) why should a predicativist believe that
the fact that some set of order type $\alpha$ is well-founded renders proof
trees of height $\alpha$ valid, and (2) granting that he can draw this
inference for any particular $\alpha$, why does he fail to grasp that it is
valid in general? The inference superficially seems reasonable
because it is classically valid, but it is hard to imagine what its
predicative justification could be. It is even harder to believe that
a predicativist could recognize its validity in each instance but not
as a general rule.

Working only in Peano arithmetic, a predicativist should be able to draw
conclusions about what is provable in some infinitary system using proof
trees of various heights. But in order to infer which proof trees are actually
valid, he needs some new principle going beyond Peano arithmetic. Augmenting
PA with an axiom schema which expresses the principle ``if a well-founded
proof tree proves $\Ac$, then $\Ac$'' in some form would yield a system
which proves the well-foundedness of a notation for \Gm{}. Expressing the
principle as a deduction rule schema rather than as an axiom schema would
yield a system which proves the well-foundedness of notations for all ordinals
less than \Gm{} but not \Gm{} itself, but we would need to explain why the
deduction rules are valid while the corresponding axioms are not (and we would
still be able to get beyond \Gm{} by a single act of reflection). One is
virtually forced to assert that whenever a predicativist proves that a set
is well-founded he is then able to infer the validity of proof trees of that
height via an unformalizable leap of intuition, but I see no reasonable basis
for such a claim.

\subsection{The linked systems $P$ and $\exists/P$}\label{sect1f}

$P + \exists/P$ can be criticized in three different ways.
\medskip

\noindent {\bf (a)} {\it Obscure formulation.} The central feature
of $P + \exists/P$, its division into two distinct but interacting
formal systems, is so unusual that it would seem to call for an
especially careful account of the underlying motivation. Although
\cite{Fef3} contains a substantial amount of prefatory material,
there is no explicit discussion of this seemingly crucial point.
One gets a vague sense that part of the motivation is to allow use
of second order quantifiers only during brief excursions into the
``auxiliary'' system $\exists/P$ as a sort of next-best alternative
to prohibiting them altogether, but nothing is said about why exactly
this degree of usage is deemed acceptable. This makes it difficult to
evaluate $P + \exists/P$, since one is left with the basic question
of how we are supposed to regard the predicative meaning and reliability
of statements proven in $P$ as opposed to those proven in $\exists/P$.

There apparently is some basic distinction to be made between the
conceptual content of the theorems of the two systems. I infer this
from the requirement both in the description of allowed formulas of
$\exists/P$ (\cite{Fef3}, p.\ 76) and in the rules IV and V
(\cite{Fef3}, p.\ 78) that
at least part of the premise specifically be proven in $P$. For
instance, the functional defining axioms (IV) allow the introduction
of a functional symbol provided existence of the functional has been
proven in $\exists/P$ and its uniqueness has been proven in $P$.
Existence can only be proven in $\exists/P$ since $P$ lacks the
necessary quantifiers, but no reason is given why uniqueness must
be proven in $P$. Would a proof of uniqueness in $\exists/P$ be
unreliable in some way? If so, why should we trust other theorems
of this system? {\it Why are we able to justify introducing a
functional symbol when the functional's uniqueness has been proven
in $P$, but not when its uniqueness has been proven in $\exists/P$?}

The question is significant because an identical point can be made in
the two other cases, and if they were all broadened to include premises
proven in $\exists/P$ then the system $P$ would become superfluous: all
reasoning could take place in $\exists/P$. Agreeing that $P$ is indeed
dispensable is not acceptable, since this would obviate the need for
the functional defining axioms altogether and thereby void Feferman's
justification for not allowing a predicativist to reflect on the validity
of $P + \exists/P$ (see \S \ref{sect1f} (c)).
\medskip

\noindent {\bf (b)} {\it Too strong.}
It is also unclear how to reconcile the proposed formalism surrounding
second order existential quantification with the motivating idea that

\begin{quotation}
\noindent we have {\it partial understanding of 2nd order existential
quantification}, for example when a function or predicate satisfying
an elementary condition is shown to exist by means of an explicit
definition. Some reasoning involving this partial understanding
may then be utilized, though 2nd order quantifiers are not to be
admitted as logical operators in general. (\cite{Fef3}, p.\ 71;
italics in original)
\end{quotation}

\noindent For example, this intuition seems somewhat incompatible with the
use of {\it negated} second order existential quantifiers, which are allowed
in $\exists/P$. Even more problematic is the proof of transfinite recursion
over well-ordered sets (\cite{Fef3}, pp.\ 82-83), which conflicts rather
severely with any understanding of second order existence in terms of
``explicit definition''. The offending aspect of this proof is
its use of predicate substitution with a $\Sigma^1_1$ formula, which
is hard to reconcile with the idea that only ``some'' reasoning about
``partially understood'' second order quantifiers is available.${}^6$
General freedom to replace set variables with $\Sigma^1_1$ formulas seems
to imply a {\it complete} ability to reason abstractly about second order
existence.

Feferman mentions the prima facie impredicative nature of his predicate
substitution rule (rule V) and responds that

\begin{quotation}
\noindent By way of justification for the schema V it may be argued that
the (predicative) provability of $\Bc(X)$ establishes its validity
also for properties whose meaning is not understood, just as one
may reason logically with expressions whose meaning is not fully
known or which could even be meaningless. (\cite{Fef3}, p.\ 92)
\end{quotation}

\noindent But this line of argument would equally well justify full
comprehension. Indeed, for any formula $\Ac(n)$ the predicatively
valid statement $(\exists Y)(n \in Y \leftrightarrow n \in X)$
yields $(\exists Y)(n \in Y \leftrightarrow \Ac(n))$ by predicate
substitution. Even if we restrict ourselves to $\Sigma^1_1$ formulas
$\Ac$, we could still infer $\Pi^1_1$ comprehension.
So the idea that ``the predicative provability
of $\Bc(X)$ establishes its validity also for properties whose
meaning is not understood'' is clearly not acceptable as a general
principle as it stands.
\medskip

\noindent {\bf (c)} {\it Too weak.}
Now consider the general objection of \S \ref{sect1c}. Feferman
addresses it in the following way:

\begin{quotation}
\noindent $\ldots$ this is not a good argument because the
functional defining axioms are only given by a generation
procedure and the predicative acceptability of these axioms
is only supposed to be recognized at the stages of their
generation. To talk globally about the correctness of $P$ we
have to understand globally the meaning of all functional
symbols in $P$; there is no stage in the generation process
at which this is available (\cite{Fef3}, p.\ 92).
\end{quotation}

\noindent The point here is that $P + \exists/P$ contains a rule which
allows one to introduce a symbol for a functional $\vec{\alpha} \mapsto
\beta$ once a unique $\beta$ satisfying some formula $\Ac(\vec{\alpha}, \beta)$
(with $\vec{\alpha} = (\alpha_1, \ldots, \alpha_n)$, and all free variables
in $\Ac$ shown) has been proven to exist for any $\vec{\alpha}$. The
$\alpha_i$ and $\beta$ are predicate variables.

I do not see how the fact that new symbols can be introduced could
in itself prevent anyone from grasping the overall validity of the
system. Surely, ``to talk globally about the correctness of $P$''
we need only to accept the validity of the functional generating
procedure, not necessarily ``to understand globally the meaning
of all functional symbols'' beforehand.${}^7$

The more substantial question is whether the validity of the functional
defining axioms of $P$ might only be recognized in stages. Now it may
be possible to imagine a set of beliefs which would lead one to accept
the functional defining axioms only at the stages of their generation:
perhaps someone could, by brute intuition, accept the validity of a
specific functional definition after having grasped an explicit
construction of the functional being defined, yet not be able to
reason about functional existence in general terms. This seems
something like the standpoint of ``immediate predicativism'' discussed
on pp.\ 73 and 91 of \cite{Fef3}. The problem is that it tends to conflict
with rule VII (relative explicit definition) and axiom VIII (unification)
(\cite{Fef3}, p.\ 78) of $\exists/P$, both of which do presume an ability
to reason abstractly about second order existence (not to mention rule V,
predicate substitution). Thus, Feferman's argument belies his premise that
a predicativist is capable of accepting rule VII and axiom VIII.${}^8$

\subsection{The system \RefP{}}\label{sect1g}

The ${\rm Ref}^*$ construction applies to any schematic formal
theory, but the case of interest for us is schematic Peano arithmetic
\PAP{}. This is formulated in the language $L$ of first
order arithmetic augmented by a single predicate symbol $P$.
The axioms are the usual axioms of Peano arithmetic with the
induction schema replaced by the single axiom
$$P(0) \wedge (\forall n)(P(n) \rightarrow P(n'))
\rightarrow (\forall n)\, P(n),$$
and there is an additional deduction rule schema allowing substitution of
arbitrary formulas for $P$. Now if ${\rm S}(P)$ is any schematic theory
then \RefS{} is a theory in the language of ${\rm S}(P)$ augmented by
two predicate variables $T$ and $F$ whose axioms are the axioms of
${\rm S}(P)$ together with ``self-truth'' axioms governing the
partial truth and falsehood predicates $T$ and $F$, and with a
substitution rule which allows the substitution of formulas possibly
involving $T$ and $F$ for $P$.
\medskip

\noindent {\bf (a)} {\it Too strong.}
In \S \ref{sect2e} and \S \ref{sect2f} I will discuss the prima facie
impredicativity both of self-applicative truth predicates and of schematic
predicate variables. Leaving those issues aside for now, the first point
to make here is that the key axiom which distinguishes \RefP{} from a
much weaker system ${\rm Ref}({\rm PA}(P))$, axiom 3.2.1 (i)${}^{(P)}$
(\cite{Fef4}, p.\ 19), has no obvious intuitive meaning. The reason for
using schematic formulas, as opposed to ordinary second order formulas
involving set variables, is that they are supposed
to allow one to fully express principles such as induction without assuming
any comprehension axioms (\cite{Fef4}, p.\ 8). This means that we
interpret a statement involving a schematic predicate symbol $P$
not as making an assertion about a fixed arbitrary set, but
rather as a sort of meta-assertion which makes an open-ended claim that
the statement will be true in any intelligible substitution instance.
However, the truth claim
$$T(\ulcorner P(\overline{n})\urcorner) \leftrightarrow P(n),$$
a special case of 3.2.1 (i)${}^{(P)}$, cannot be given the latter
interpretation since the number $\ulcorner P(\overline{n})\urcorner$
does not change when a substitution is made for $P$ in this formula. If
we interpret $P$ in a way that is compatible with 3.2.1 (i)${}^{(P)}$,
i.e., as indicating membership in a fixed set, then the substitution rule
${\it P-Subst:}\, L(P)/L(P,T,F)$ (\cite{Fef4}, Definition 3.3.2 (iii))
can only be justified by an impredicative comprehension principle
(cf.\ \cite{Fef4}, p.\ 8).

Feferman characterizes axiom 3.2.1 (i)$^{(P)}$ as ``relativizing $T$
and $F$ to $P$'' (\cite{Fef4}, p.\ 19). I am not sure what this means,
but the axiom clearly is not valid on arbitrary
substitutions for $P$, yet one draws consequences from it to
which one does apply a substitution rule (and this is crucial
for the proof that $(\Pi^0_1-{\rm CA})_{<\Gamma_0} \leq$ \RefP{}).
In fairness, I should point out that this problem is noted in
(\cite{Fef4}, \S 6.1.3 (i)), with the comment that ``a
fall-back line of defense could be that this substitution
accords with ordinary informal reasoning. However, this seems
to me to be the weakest point of the case for reflective
closure having fundamental significance.''

I would argue that the above difficulty not only invalidates the
idea that \RefP{} models predicative reasoning, it shows that the
${\rm Ref}^*$ construction indeed has no fundamental significance.
There is no way to interpret the $P$ symbol that simultaneously makes
sense of the axiom 3.2.1 (i)${}^{(P)}$ and the substitution rule
${\it P-Subst:}\, L(P)/L(P,T,F)$.
\medskip

\noindent {\bf (b)} {\it Too weak.} The ${\rm Ref}^*$ construction is
described in \cite{Fef4} as a ``closure'' operation and the question
of its significance is discussed in terms of Kripke's theory of
grounded truth outlined in \cite{Krip}. A casual reading of \S 6
of \cite{Fef4} might leave the impression that the statements $\Ac$ such
that \RefP{} proves $T(\ulcorner{\Ac}\urcorner)$ are supposed to be
precisely the grounded true statements of the language $L(P, T, F)$. But
this cannot be right because these statements are recursively enumerable,
so that one can write a formula $(\forall n)\,T(\{\overline{r}\}(n))$ which
asserts precisely their truth. This formula is grounded (in any
reasonable sense) and true but the assertion of its truth is not
a theorem of the system.

The more careful formulation that the self-truth axioms ``correspond
directly to the informal notion of grounded truth and falsity''
(\cite{Fef4}, p.\ 42) is well-taken, but we must not confuse this
with the claim that the self-truth axioms {\it capture} the informal
notion. Their failure to do so can be traced to axiom (vi)
(\cite{Fef4}, Definition 3.2.1). The first conjunct of this axiom,
for example, asserts that if $\Ac(\overline{n})$ is true for every $n$ then
$(\forall n)\,\Ac(n)$ is true. But this does not fully capture the
informal idea that ``the truth of $\Ac(n)$ for all $n$ implies the
truth of $(\forall n)\,\Ac(n)$'' in the sense that there exist
formulas $\Ac(n)$ which can be proven true for each numerical value
of $n$ but such that there is no proof in \RefP{} that for all $n$,
$\Ac(\overline{n})$ is true --- in particular, $\Ac(n) \equiv T(\{\overline{r}\}(n))$
where $\{r\}$ enumerates all $a$ such that \RefP{} proves
$T(\overline{a})$.

Now consider the claim that in general \RefS{} encapsulates what
one ``ought to accept'' given that one has accepted ${\rm S}(P)$
(\cite{Fef4}, p.\ 2).
This has an air of paradox since one has to ask whether the claim
itself is something that anyone ought to accept. However, that
point is not crucial to the question of what a predicativist
can prove since it need not attach to the specific assertion that
\RefP{} encapsulates what one ought to accept given that one has
accepted \PAP{}. We may suppose that the predicativist
does not realize (and indeed, {\it ought not accept}) that his
commitment to Peano arithmetic obliges him to accept every theorem
of \RefP{}, although this is in fact the case. This leads us back
to the question posed in \S \ref{sect1c}. Evidently we are
dealing with a claim that predicativists can affirm each theorem
of \RefP{} individually but cannot accept this system globally.

This point is not explicitly addressed in \cite{Fef4}, but the
informal notion of ``partial truth'' has the flavor of a forever
incomplete process and might seem like it could support such a claim.
For example, it is suggested in \S 6.1.1 of \cite{Fef4} that
the passage from ${\rm S}(P)$ to ${\rm Ref}({\rm S}(P))$ should not
be iterated because this would ``vitiate the informal idea behind
the use of partial predicates of truth and falsity.'' A possibly more
straightforward question which avoids the issue of using multiple
partial truth predicates is whether one could justify augmenting
\RefP{} by the single statement $(\forall n)\,T(\{\overline{r}\}(n))$
described above.

Surely a predicativist {\it can} justify adding this statement
if he is able to generally recognize that every statement proven true by
\RefP{} is indeed true. Given that \RefP{} is finitely axiomatized and
that the predicativist is presumed to accept each theorem of \RefP{}
individually, it is unclear how this could be plausibly denied. Indeed,
axiom (vi) clearly affirms that the predicativist {\it is} able to reason
about the collective truth of an infinite set of statements themselves
involving assertions of truth and falsehood.

In \S 6.1.3 (ii) of \cite{Fef4} Feferman considers the question ``have we
accepted too little?'' in terms of logically provable statements, e.g.\ of
the form $\Ac \vee \neg \Ac$, whose truth is not provable because $\Ac$ is
not grounded. This leads into a brief discussion of the relative merits of
Kripke's minimal fixed point approach versus van Fraassen's more liberal
``supervaluation'' approach to self-applicative truth. But this discussion
is misleading because \RefP{} does not even prove the truth of every
statement in Kripke's minimal fixed point; in particular, if this needs
repeating, it does not prove the statement
$(\forall n)\, T(\{\overline{r}\}(n))$.
This formula is not logically provable but it is grounded, and it seems a
rather stronger candidate for a statement that ``ought'' to be accepted
as true.

\subsection{The system \UNFA{}}\label{sect1h}

Distinct, not obviously equivalent, versions of \UNFA{} are presented in
\cite{Fef6} and \cite{FefSt}. I give priority to the later version in 
\cite{FefSt}.

\medskip

\noindent {\bf (a)} {\it Too weak.} Like ${\rm Ref}^*$, $\Uc$ is
presented in \cite{FefSt} as a general construction (``unfolding'')
which can be applied to any schematic formal system ${\rm S}(P)$.
As usual, granting that acceptance of ${\rm S}(P)$ justifies
acceptance of every theorem of $\Uc({\rm S}(P))$, we can ask why
it fails to justify accepting a formalized $\omega$-rule schema
referring to theorems of $\Uc({\rm S}(P))$. This question is not
addressed in either \cite{Fef6} or \cite{FefSt}; the closest
I can find to an answer is the following passage in \cite{Fef6}:

\begin{quotation}
\noindent [W]e may expect the language and theorems of the
unfolding of (an effectively given system) ${\rm S}$ to be
effectively enumerable, but we should not expect to be able
to decide which operations introduced by implicit (e.g.\
recursive fixed-point) definitions are well defined for all
arguments, even though it may be just those with which we
wish to be concerned in the end. This echoes G\"odel's
picture of the process of obtaining new axioms which are
``just as evident and justified'' as those with which we
started $\ldots$ for which we cannot say in advance exactly
what those will be, though we can describe fully the means
by which they are to be obtained. (\cite{Fef6}, p.\ 10)
\end{quotation}

\noindent Here reference is made to the fact that $\Uc$ uses
partial operations, which is seen as having
fundamental significance. It is true that the question of
which partial operations of \UNFA{} are total is (unsurprisingly)
not decidable, though this in itself seems a questionable basis
for forbidding us from proceeding beyond \UNFA{} when it did not
prevent us from formulating this system in the first place or
from working within it.

If ${\rm S}(P)$ involves no basic objects of type 2 (as is the case for
${\rm NFA}$) then an argument could be made that applying the $\Uc$
construction twice is conceptually different from applying it once in
that $\Uc({\rm S}(P))$ does employ higher type objects and thus the
original system ${\rm S}(P)$ can possibly be seen as being ``concrete''
in a way that $\Uc({\rm S}(P))$ is not. However, this should not prevent
one from accepting a formalized $\omega$-rule schema applied to \UNFA{},
which would seem to require only that one accept \UNFA{} is sound.${}^9$
\medskip

\noindent {\bf (b)} {\it Way too strong.} \UNFA{} is actually
flatly impredicative in two distinct ways. First, the $\Uc$
construction suffers from the same nonsensical treatment of
schematic predicates as ${\rm Ref}^*$. Here the offending axiom
is Ax 7 (\cite{FefSt}, p.\ 82), which does not make sense if
$P$ is understood as a schematic predicate. It is valid if
we regard $P$ as indicating membership in a fixed set,
but then, just as for \RefP{}, use of the substitution
rule (Subst) (\cite{FefSt}, p.\ 82) would have to presume
impredicative comprehension.

The really striking impredicativity of \UNFA{}, however, is its use
of a least fixed point operator, which apparently informally assumes
the legitimacy of generalized inductive definitions in the sense
of \cite{BFPS}. This not only vitiates any claim of \UNFA{} to model
predicative reasoning, it more broadly undermines the idea that \UNFA{}
has any fundamental philosophical significance, since it would seem
that anyone who accepts the $\Uc$ construction and Peano arithmetic
ought to at least accept ${\rm ID}_1$ \cite{BFPS}, which is far
stronger than \UNFA{}.${}^{10}$

\subsection{Summary of the critique}\label{sect1i}

At the beginning of this section I made strong claims about the
weakness of the case for the \Gm{} thesis. Were they borne out?

First, I stated that each of the formal systems of \S \ref{sect1a}
is motivated by informal principles which actually justify a stronger
system that proves the well-foundedness of an ordered set that is
isomorphic to \Gm{}. In the case of $\Sigma$, $H$, and $R$, the informal
principle is ``$a$ is an ordinal notation implies $S_a$ is sound'', which
is needed to justify ($*$) but in fact justifies ($**$) (see \S
\ref{sect1d}). In $H^+$, $R^+$, and $RA^*$ the principle is ``$a$ is
an ordinal notation implies proof trees of height $a$ are sound''.
In $P + \exists/P$ we accept
that it is legitimate to substitute arbitrary predicates for set
variables, which justifies full comprehension. \RefP{} assumes an
informal grasp of a self-applicative concept of truth, which
justifies the inference of a statement that asserts the truth of every
theorem proven true by \RefP{}. \UNFA{} informally assumes the legitimacy
of generalized inductive definitions, which actually justifies ${\rm ID}_1$.

Second, I stated that the responses on record to the objection in
\S \ref{sect1c} are brief, vague, and unpersuasive. The only such
responses of which I am aware are Kreisel's answers in \cite{Kre1},
\cite{Kre5}, and \cite{Kre6} (see \S \ref{sect1d} (a), (b), (c)) and
Feferman's answer about $P + \exists/P$ in \cite{Fef3} (see \S \ref{sect1f}
(c)). In \cite{Kre1} and \cite{Fef3} the response is barely more than a
flat assertion with no real explanation given; in \cite{Kre5} it is a
cryptic passage whose most reasonable interpretation is clearly
self-defeating; and in \cite{Kre6} it is merely an implausible
speculation. With regard to \RefP{} and \UNFA{}, as far as I am aware
the objection has not even been discussed in the literature, except
tangentially by an argument in \cite{Fef4} that the ${\rm Ref}^*$
construction should not be iterated.

Finally, I said that each of the formal systems of \S \ref{sect1a}
is manifestly impredicative in some way. The most blatant example of
this is the least fixed point operator in \UNFA{}, but in all three
of $P + \exists/P$, \RefP{}, and \UNFA{} there is a basic impredicativity
involving the ability to substitute possibly meaningless formulas
for free set variables. In \RefP{} and \UNFA{} this is hidden by employing
a substitution rule involving a ``schematic'' predicate
symbol, but elsewhere treating this predicate symbol in a way that only
makes sense if it is thought of as a classical predicate indicating
membership in a fixed set.${}^{11}$

As I mentioned in \S \ref{sect1b}, the reason one needs a substitution
rule is because one wants to convert statements of transfinite induction
into statements of transfinite recursion so that one can construct iterative
hierarchies. In the autonomous systems this is accomplished by simply
postulating that a statement of transfinite induction legitimates a
transfinite application of reflection principles, which allows one to
pass to a stronger system that proves the existence of the next hierarchy.
Thus, {\it every} system uses an impredicative step to get from transfinite
induction to transfinite recursion. This is not surprising, as there is a
predicatively essential difference between induction and recursion (see the
end of \S \ref{sect2d}).

There is still a question as to why so many different attempted
fomalizations of predicative reasoning happen to have proof-theoretic
ordinal \Gm{}. One answer is that they all employ essentially the same
well-ordering proof and to a substantial extent appear to have been {\it built
around} different versions of this proof. (This is especially evident
in the case of $P + \exists/P$.) Another possible answer is that
the earlier systems all access the same ordinals because they embody
the same fallacies revolving around the idea of autonomous generation
of ordinals, while the later systems were formulated against the
background of the earlier systems which had already seemed to attain
the correct answer. This could have made it difficult to free oneself
from a conclusion that had already been formed and seemed well-justified.
However, a properly functioning scientific community should be expected
to debate and criticize major ideas, not to passively accept them,
regardless of the stature of their author and the complexity of the
argument. That this apparently was not done in a serious way in the
present case suggests that the community as a whole did not function
as it should have.

\section{An analysis of predicative provability}\label{sect2}

As I discussed in \S \ref{sect1b} and \S \ref{sect1i}, all of the formal
systems of \S \ref{sect1a} employ impredicative methods
in order to pass from transfinite induction to transfinite recursion.
This presents a basic obstacle to obtaining predicative ordinals by
means of the general technique employed by those systems. Our goal in
this section is to develop new methods of producing predicative
well-ordering proofs.

Considering the variety of impredicative features that have appeared
in previous attempts to model predicative reasoning (and there are two
other major ones besides those I discussed in Section \ref{sect1};
see \S \ref{sect2e} and \S \ref{sect2f}), it seems fair to say that
not enough attention has been paid to the basic conceptual content
of predicativism. Therefore our discussion will incorporate a general
conceptual analysis of predicativist principles.

As many commentators have noted, the vicious-circle principle --- generally
taken as the defining principle of predicativism --- does not in itself
constitute a well-defined foundational stance, as it is compatible with a
variety of attitudes about which principles of set construction should be
accepted as basic. Indeed, the version of predicativism under consideration
here is nowadays often referred to as ``predicativism given the natural
numbers'', a phrasing which I find unfortunate, as it gives no indication
as to why one should take the natural numbers as basic, as opposed to any
other set. I take the essential basis for this view to be a conception
of mathematical reality according to which sets have no independent prior
existence
but must be constructed, together with the idea that infinite constructions
are legitimate, but only if one is able to conceive them in a completely
precise and explicit way. The special role that the natural numbers play
in this account arises from the fact that ``we have a complete and clear
mental survey of all the objects being considered,
together with the basic [order] interrelationships between them''
(\cite{Fef3}, p.\ 70). Elsewhere I term this stance
{\it mathematical conceptualism} and argue in support of it \cite{W1}.
For our purposes here we may encapsulate it in three basic principles:

\begin{quotation}
\noindent (i) the mathematical universe is a variable entity
that can always be enlarged
\medskip

\noindent (ii) every set must be constructible from logically
prior sets
\medskip

\noindent (iii) constructions of length $\omega$ are legitimate.
\end{quotation}

I take it as granted that these principles consistute a coherent
foundational stance and I am not concerned here with trying to
justify them. However, a brief explanation is in order. The idea
of (i) is that there is no well-defined complete universe of sets
because any particular collection of sets can itself be identified as
a set that does not belong to the collection. Any currently available
partial universe can always be extended, and this extension process
can be iterated; since we accept constructions of length $\omega$
it can even be iterated transfinitely. But the general concept of
extension cannot be fully formalized since one can always go one
step beyond any given partial formalization.

Assertion (ii) is an informal version of the vicious circle
principle. Although that principle is notoriously difficult to
formulate precisely (in particular, we do not attempt to define
``logically prior''), the underlying intuition seems clear enough
to be used in most cases in evaluating whether a proposed formal
system is predicatively acceptable.

The idea behind (iii) is that we have an intuitively clear
conception of what it would be like to carry out a process (construction,
computation, proof) of length $\omega$ and therefore the results of
such procedures are legitimate objects of study. By iteration we can
accept processes of length $\omega^2$, $\omega^\omega$, etc., but
for now we leave open just what the ``etc.''\ entails. As with (i)
we do not expect to be able to fully formalize exactly how far we
can go.

\subsection{The power set of $\omega$}\label{sect2a}

Unlike naive set theory, predicativism obviously does not support the
principle ``for any property $\Ac$ of sets, $\{X: \Ac(X)\}$ exists''.
Indeed, the intuitive appeal of this general comprehension principle
seems to rest on an implicit belief in the existence of a well-defined
complete universe of sets. If there were such a universe $V$, then
for any definite property $\Ac$ one might have some idea of forming
$\{X: \Ac(X)\}$ by extracting from $V$ just those sets satisfying
$\Ac$. But if we reject the existence of such a universe then this
idea fails, and in fact for many properties $\Ac$ (e.g., ``$0 \in X$'')
we clearly {\it cannot} imagine any way to form the set of all $X$'s
which satisfy $\Ac$. It is important to understand that this can be
so even if $\Ac$ is {\it definite} in the sense that $\Ac(X)$ is
recognized to have a well-defined truth value for any conceivable $X$.

Thus, we do not accept a set as legitimate if it can only be defined
``from above'' in the form $\{X: \Ac(X)\}$. We do accept sets which
can be defined by restricted comprehension (i.e., are of the form
$\{X \in Y: \Ac(X)\}$) relative to a set which has already been accepted,
provided the property $\Ac$ is definite in the sense just indicated,${}^{12}$
and we also accept sets which can be constructed ``from below''.
Principle (iii) gives us a powerful ability to build up countable sets
from below, but we would not expect that we can predicatively reach
any classical uncountable set in this manner. Indeed, on the basis of
the principles enumerated above it is reasonable to assume that every
predicatively acceptable infinite set should be not only classically
but predicatively countable (i.e., predicatively known to be in bijection
with $\omega$).

For our purposes here we need not insist on such an ``axiom of
countability'', but we have to accept that, at the very least, we cannot
assume that the power set of $\omega$ is a predicatively legitimate set.
As explained above, this does {\it not} contradict the definiteness of
the property ``$X$ is a subset of $\omega$'', and indeed the latter
could be justified by an appeal to (iii). Given any set $X$, an informal
``computation'' of length $\omega +1$ could verify or falsify the claim
that $X \subseteq \omega$: for each $n \in \omega$ check whether $n$ belongs
to $X$; if so, remove it; at step $\omega$ check whether any elements
remain. Thus, the property of being a subset of $\omega$ is
predicatively definite.${}^{13}$

\subsection{Predicatively valid logic}\label{sect2b}

Classical logic is unsuited to reasoning about a variable
universe. Since the general extension process by which new sets are
recognized cannot be completely formalized, the universe is inherently
ill-defined and so we do not expect every assertion about sets to have
a well-defined truth value. Rather, we should regard the family of true
statements as a variable entity which is always capable of
enlargement, much like the mathematical universe itself. {\it This
makes intuitionistic logic the appropriate tool for general
predicative reasoning.}

Of course, this is not to say that the predicative notion of truth
can be identified with intuitionistic truth. Predicatively there is
no reason to believe (and undoubtedly good reason not to believe) that
every true statement can be proven by a finite argument. Conversely,
for reasons I do not understand,
it seems that most intuitionists accept impredicative constructions.
Nevertheless, I maintain that the logical apparatus of intuitionism
is exactly suitable for predicativism. To say that the law of the
excluded middle always holds is just to say that {\it every} formula
is definite in the sense of \S \ref{sect2a}. Predicatively the
definiteness of any statement that quantifies over subsets of
$\omega$ is initially suspect, so it is highly implausible that a
predicativist could be led to accept even that all formulas of second
order arithmetic are definite.

For a specific example of the presumable failure of the law of the
excluded middle, notice that well-ordering assertions can apparently
fail to have a well-defined truth value because the inherent ambiguity
of the mathematical universe could lead to uncertainty about whether or
not a given totally ordered set has a proper progressive subset. If no
such subset is currently available, indefiniteness about whether such a
subset will appear in some future enriched universe could be a reasonable
consequence of the fact that {\it we do not know how new sets might
arise}. An even sharper example is given by the set
$S = \{n: \Ac_n$ is true$\}$ where $(\Ac_n)$ is some recursive
enumeration of the sentences of second order arithmetic. The set $S$ is
obviously impredicative since it is defined in terms of quantification over
the power set of $\omega$, but if we accepted $(\forall n)(\Ac_n \vee
\neg \Ac_n)$ then the restricted comprehension principle mentioned in
\S \ref{sect2a} would allow us to form $S$. This shows that
$\Ac \vee \neg \Ac$ must not be assumed to hold in every case.

On the other hand, we do regard statements relativized to any
well-defined partial universe as definite, so that any such statement
should have a well-defined truth value. For example, in the setting of
second order arithmetic, principle (iii) should at least assure us that any
arithmetical statement is definitely true or false since we can imagine
checking it mechanically. This is so even if the statement contains set
variables as parameters, since for any particular $n \in \omega$ and $X
\subseteq \omega$ the atomic formula ``$n \in X$'' has a definite truth value.
Thus, at the level of arithmetical statements our logic is classical.

Similar considerations were discussed in \cite{CoPa}, leading to the
suggestion that predicativists can adopt the {\it numerical omniscience
schema}
$$(\forall n)\left(\Ac(n)\vee \neg \Ac(n)\right) \rightarrow
\left[(\forall n)\,\Ac(n) \vee (\exists n)\, \neg \Ac(n)\right]$$
(where here $\Ac$ is any formula of second order arithmetic and $n$ is
a number variable). Together with the assumption $\Ac \vee \neg\Ac$
for every atomic formula $\Ac$, this implies the law of the excluded
middle for every arithmetical formula.

A word about terminology. If we do not assume the law of the excluded
middle then we may have to consider assertions whose sense is understood
but which are not known to have a definite truth value. To keep this
distinction clear I will say an assertion is {\it meaningful} if it has
a definite truth value and {\it intelligible} if its sense is understood.
Thus, every meaningful assertion must be intelligible and every intelligible
assertion is potentially meaningful.

\subsection{Second order quantification}\label{sect2c}

First order (numerical) quantification is unproblematic by principle
(iii). The legitimacy of second order quantification is less clear
since we do not regard the power set of $\omega$ as a well-defined
entity over which set variables could be imagined ranging. This has
been a recurrent concern in the literature on predicativity. For
instance, it was cited as motivation for the strong restrictions
on second order quantification in \cite{Fef3}.

To what extent, if any, are second order quantifiers
acceptable? First, because the concept ``set of numbers'' is
definite (\S \ref{sect2a}), we should at least be able to make some limited
constructive sense of existential quantification. There are situations
in which we can recognize that we are (in principle) able to construct a
set of numbers with some property, and this should license some use of
second order existential quantifiers. This was also the position taken
in \cite{Fef3}.

In addition, we do seem to be able to predicatively accept some statements
as being true of any set of numbers. Despite the unfixed nature of
the mathematical
universe, we can still affirm general assertions like $0 \in X$ $\vee$
$0 \not\in X$ as holding for any conceivable $X \subseteq \omega$. Not only
is this statement true for all currently available sets, it must remain
true in any future universe. We can be sure that we will never
come across a set for which the assertion fails because its truth is
inherent in the concept ``set of numbers''. As another example, given any
$X \subseteq \omega$, principle (iii) should justify asserting the
(constructive) existence of its complement. Thus, we ought to be able
to somehow express that for every $X$ there is a $Y$ such that $n \in Y
\leftrightarrow n \not\in X$. Finally, the principle of induction in the
form $0 \in X \wedge (\forall n)(n \in X \rightarrow n' \in X) \rightarrow
(\forall n)(n \in X)$ is recognizably true for any $X \subseteq \omega$.
Given that we accept processes of length $\omega$, we can be certain that
any set which satisfies the induction premise must contain every number
since we can imagine verifying this conclusion mechanically. Again, this
must hold not only for all currently available sets but for all sets
in any conceivable future universe.

In \cite{Fef9}, following Russell, a distinction is drawn between
the concepts ``for all'' (ranging over a well-defined collection)
and ``for any'' (ranging over a ``potential totality''). I find
this distinction helpful, but in the present setting I do not
accept Russell's suggestion, adopted in \cite{Fef3}, that the
``for any'' intuition is captured by using free set variables.
Consider the following example. We have already agreed that
predicativists can acknowledge that any subset $X$ of $\omega$
has a complement $Y$. But they should then also agree that $Y$
has properties like: for any $Z$, $Z\subseteq Y$ if and only
if $Z \cap X = \emptyset$. Indeed, given any $X$ and $Z$ we can
imagine constructing $Y$ (using (iii)) and then verifying the
relation between $X$, $Y$, and $Z$ (again using (iii), specifically
a version of the numerical omniscience schema, together with definiteness
of the assertions $Z \subseteq Y$ and $Z \cap X = \emptyset$). Since the
construction of $Y$ did not depend on $Z$ this means that we can affirm
the statement
$$(\forall X)(\exists Y)(\forall Z)(Z \subseteq Y \leftrightarrow
Z \cap X = \emptyset)$$
under the interpretation $\forall$ = ``for any'' and $\exists$ =
``there can be constructed a''. This shows that alternating second order
quantifiers can make predicative sense. Moreover, the idea {\it cannot}
be expressed without using at least one universal quantifier, which
shows that Russell's free variable suggestion is inadequate here.

An alternative possibility is to allow use of both universal and existential
second order quantifiers and to reason using an intuitionistic predicate
calculus. Given the conception of predicativism developed above and
the interpretation of second order quantifiers just indicated, this
logical apparatus appears perfectly acceptable. {\it Intuitionistic
logic legitimates the predicative use of set quantifiers.}

I have been careful to restrict this discussion to sets of numbers.
Using quantifiers to range over {\it all sets}, with no qualification,
is less tenable because the general concept of a set may be
predicatively unclear. To a predicativist, the assertion that some object
is a set may not only be indefinite (in the sense of \S \ref{sect2a}) but
unintelligible (in the sense of \S \ref{sect2b}). In fact, I believe
that quantification over all sets is not predicatively valid for just
this reason. However, I would still reject Russell's free variable
convention here; rather, I would conclude that predicativists simply
cannot make general statements about all sets.${}^{14}$

(There might be ways around this difficulty. For instance, if one
interprets sets in terms of well-founded trees coded by subsets of
$\omega$, cf.\ \S VII.3 of \cite{Sim}, then the concept may become
predicatively intelligible.)

\subsection{Predicative well-ordering}\label{sect2d}

In the previous section I justified a second order induction
statement using principle (iii). For which formulas $\Ac$ of second
order arithmetic would a similar argument lead us to accept $\Ac(0) \wedge
(\forall n)(\Ac(n) \rightarrow \Ac(n')) \rightarrow (\forall n)\, \Ac(n)$?

If it contains set variables, the formula $\Ac(n)$ might not have a
definite truth value (\S \ref{sect2b}). However, once we have proven
$\Ac(0)$ we must at least agree that this instance is definitely true.
If, moreover, we have proven $(\forall n)(\Ac(n) \rightarrow \Ac(n'))$
then we can be successively brought to the same conclusion about $\Ac(1)$,
$\Ac(2)$, etc., and recognizing this, we should therefore accept
$(\forall n)\, \Ac(n)$ as true. Regarding the family of true statements
as a variable entity always capable of enlargement, this shows that
predicativists should accept induction for {\it every} formula $\Ac(n)$.

This may need further explanation in light of my insistence in \S
\ref{sect1f} (b) that it is generally not valid to substitute arbitrary,
possibly meaningless, formulas for set variables. I stand on this
assertion: for example,
$(\forall n)(n \in X$ $\vee$ $n \not\in X)$ is predicatively true
but $(\forall n)(\Ac(n) \vee \neg \Ac(n))$ is presumably not if, e.g.,
$\Ac(n)$ asserts that $n$ is a Church-Kleene ordinal notation. However,
this does not entail that possibly meaningless formulas can never
appear in true statements. A predicativist should accept complete
induction (provided he is using intuitionistic logic)
since he can generally recognize that the truth of the premise of
any induction statement would entail the truth of its conclusion
even if the latter was not initially known to be meaningful.

Next let us consider the extent to which predicativists can understand
the general concept of a well-ordered set. It is sometimes said that the
well-ordering concept is not available to predicativists because it
involves quantification over power sets. On the other hand, it seems to
be generally accepted that predicativists are able to assert relatively
strong versions of the statement that $\omega$ is well-ordered.
If we agree with the conclusions of \S \ref{sect2c} then statements
of transfinite induction of the form $(\forall X)\, {\rm TI}(X,a)$
(transfinite induction up to $a$ on a totally ordered subset of $\omega$)
are predicatively intelligible.${}^{15}$ Here I use the abbreviations
\begin{eqnarray*}
{\rm TI}(X,a)
&\equiv&{\rm Prog}(X) \rightarrow
(\forall b \prec a)(b \in X)\\
{\rm Prog}(X)
&\equiv&(\forall b)\left[(\forall c \prec b)(c \in X)
\rightarrow b \in X\right].
\end{eqnarray*}

I argued above that complete induction on $\omega$ is predicatively valid.
Note, however, that if we know $\{b: b \prec a\}$ is well-ordered, i.e.,
we have verified $(\forall X)\, {\rm TI}(X,a)$, we cannot in general infer
${\rm TI}(\Ac,a)$ ($\equiv {\rm Prog}(\Ac) \rightarrow (\forall b \prec a)
\Ac(b)$ where ${\rm Prog}(\Ac) \equiv (\forall b)[(\forall c \prec b)
\Ac(c) \rightarrow \Ac(b)]$) for arbitrary formulas $\Ac$. The
latter schema is genuinely stronger because $(\forall X){\rm TI}(X,a)$
only asserts induction for {\it sets} that are by assumption well-defined,
whereas ${\rm TI}(\Ac,a)$ can hold if $\Ac$ is not meaningful,
and it can even be used to prove that $\Ac(b)$ {\it is} meaningful
for all $b \prec a$. It may in fact be the case that whenever
there is a predicatively valid proof of $(\forall X)\, {\rm TI}(X,a)$
there is also a proof of ${\rm TI}(\Ac,a)$ for any intelligible formula
$\Ac$. However, inferring the latter statement from the former seems to
me clearly predicatively illegitimate.

\subsection{Schematic assertions}\label{sect2e}

Even the complete induction schema does not entirely capture
a predicative understanding of induction on $\omega$ since it
only covers formulas that can be written in the language
that is currently in use. If the expressive power of the
language were strengthened in a predicatively intelligible way,
then a predicativist should accept the induction schema for all
formulas of the new language too.

This issue is addressed in \cite{Fef4} and \cite{Fef6} by a proposal
to use a ``schematic'' predicate symbol $P$ and to express the principle
of induction in a single schematic formula. Together with an informal
commitment to continue to accept all substitution instances of this
statement if the language is enriched in any intelligible way, this
does seem to fully capture a predicative understanding of induction
on $\omega$. However, it seems unlikely that a predicativist could
agree to accept such a schematic formulation because of the circularity
involved in having a formula $\Ac(P)$ which contains a schematic
predicate symbol $P$ that ranges over a class of formulas that
includes $\Ac(P)$. In fact, without some special justification
this usage seems clearly impredicative.

There should be no problem in using a schematic predicate symbol to range
over all formulas of a previously accepted language, or even a previously
accepted set of languages, as this would present no possibility of
circularity. However, because of the inherently impredicative quality
of a self-applicative predicate variable it seems to me that the general
concept ``intelligible predicate'' is itself not intelligible and that it
is therefore not possible for a predicativist to legitimately make
assertions about all intelligible predicates (cf.\ \S \ref{sect2c}).
This leads to the conclusion that {\it predicativists have an open-ended
ability to affirm induction statements on $\omega$ but are not capable of
formally expressing this fact}.

The difficulties involved with schematic predicates shed light on the
predicative unacceptability of some formal systems which superficially
have a strong predicative flavor. For example, in \cite{Kre2} and
\cite{Kre4} the possibility is raised that under intuitionistic logic
theories of generalized inductive definitions might be predicatively
valid, and this idea does have superficial appeal. However, on close
examination there is a clear circularity even in the intuitionistic
case. This is seen as follows.

Suppose we want to introduce a predicate symbol for the class defined by some
inductive definition. Classically we could define this class ``from above''
as the intersection of all classes satisfying the relevant closure condition,
but this is clearly impredicative. In the intuitionistic setting we instead
conceive of the class as an incomplete entity that can always be enlarged by
repeatedly applying the closure condition, which seems to be a predicatively
legitimate idea. The problem is in verifying the minimality property of this
class. Let $\Ac$ be any formula in the language of first order arithmetic
enriched by a predicate symbol $I_X$ which is to represent the class $X$
being defined; assuming $\Ac$ satisfies the same closure condition as $X$, we
must affirm $(\forall n)(I_X(n) \rightarrow \Ac(n))$. Now what is immediately
clear from our conception of $X$ is that this statement is progressive in the
sense that  if it holds at all previous stages in the construction of $X$
then it will still hold at the immediately following stage since $\Ac$
satisfies the same closure condition as $X$. This suggests that the
statement should be verified by a transfinite induction and we must therefore
imagine the stages in the construction of $X$ as corresponding to elements
of a well-ordered set. The difficulty then lies in specifying what we
mean by ``well-ordered''. If we had the ability to make assertions like
${\rm TI}(P,a)$ where $P$ is a schematic predicate variable, then we could
take ``well-ordered'' to mean ``supports transfinite induction for a
schematic predicate'', and we should then be able to carry out the
transfinite induction needed to prove minimality. But if the most we can
say of any totally ordered set is that it supports transfinite induction
for all formulas of a given previously accepted language, then $X$ cannot
be conceived as being built up along sets that support transfinite
induction for formulas of a language that includes $I_X$. This
would be circular because the well-ordering assertion would refer
to the class $X$ which it is being used to define. But proving the
minimality statement requires that we be able to carry out transfinite
induction for formulas of this language. Hence there is no (or at least
no obvious) way to predicatively verify minimality.

Kripke-Platek set theory also has a superficial predicative plausibility,
assuming that the point raised at the end of \S \ref{sect2c} does not
preclude any predicative treatment of arbitrary sets, but it fails for
a similar reason. Namely, the KP foundation schema is
impredicative for essentially the same reason that inductive definitions
are. For a statement of the KP axioms that is intuitionistically suitable,
see, e.g., \cite{Avi}. Their intuitionistic justification involves a
conception of an incomplete universe of sets which is built up in stages.
In order to verify any instance of the foundation schema we would
therefore need to carry out a transfinite induction with respect to
the well-ordered sets along which this universe is being constructed.
But the formula being proven by induction is a formula of the language
of KP and would implicitly make reference to the universe
being defined. Thus, in order to verify the foundation schema we would
need to build up the KP universe along sets that are known to be
well-ordered with respect to formulas which refer to that universe.
Again, this is circular and hence impredicative.${}^{16}$

\subsection{Truth theories}\label{sect2f}

Without using some kind of reflection principle I doubt that
predicativists can get beyond ordinals in the neighborhood of
$\gamma_2$ or $\gamma_3$, at most. In order to progress significantly
further we need a systematic way of iterating the process of
reflecting on the truth of a given theory to get a slightly stronger
theory. One might hope to do this using a self-applicative truth
predicate, as in \cite{Fef4}. On its face, the predicative legitimacy
of a self-applicative truth theory is problematic --- indeed, this seems
just the sort of thing that predicativist principles tend to forbid.
We can try to get around the prima facie circularity of such a theory
by regarding the truth predicate as partial and built up in stages,
giving it the flavor of a generalized inductive definition. Now
I argued in \S \ref{sect2e} that theories of generalized inductive
definitions are impredicative, but the difficulty with such theories
is their assertion of minimality axioms, which we do not require of a
truth predicate. On the contrary, the concept of belonging to an
inductively defined class does not seem predicatively objectionable
on its own; for example, according to \S \ref{sect2d} the assertion
``$n$ is a Church-Kleene ordinal notation'' is predicatively intelligible.
A parallel could also be drawn with the predicative conception of
the power set of $\omega$ as a necessarily incomplete entity that can
always be enlarged. Therefore, it seems that {\it provided intuitionistic
logic is used} self-applicative truth theories could be predicatively
justifiable. The systems of \cite{Fef4} are firmly embedded in classical
logic, but I suppose it is likely that there is an intuitionistic version
of, say, the ${\rm Ref}({\rm PA})$ construction of \cite{Fef4} that could
be accepted as predicatively legitimate. However, such a theory would
presumably have proof-theoretic ordinal only in the neighborhood of
$\gamma_2$ or $\gamma_3$. So self-applicative truth theories do not seem
a promising route to obtaining strong predicative well-ordering proofs.

Perhaps surprisingly, I find that it is possible to predicatively prove
relatively strong well-ordering assertions using hierarchies of Tarskian
(i.e., non self-applicative) truth predicates. The remainder of this paper
will develop this approach.

For the sake of readability I will begin by defining a single-step
Tarskian truth theory. Let $\Zi$ be the theory in the language of second
order arithmetic with (1) the axioms and rules of a two-sorted intuitionistic
predicate calculus and (2) the Peano axioms including induction for all
formulas of $\Zi$. We do not assume any comprehension axioms. {\it For the
remainder of the paper all theories will be assumed to have the axioms and
rules of a two-sorted intuitionistic predicate calculus and to be expressed
in the language of second order arithmetic possibly extended by a countable
family of unary relation symbols.} For each $n$ fix
a recursive bijection $\langle \cdot, \ldots, \cdot\rangle$ from $\omega^n$
to $\omega$ with corresponding recursive projections $\pi_i = \pi_i^n$,
so that $\pi_i(\langle a_1, \ldots, a_n\rangle) = a_i$. I will write
$(a_1, \ldots, a_n)$ for $(\langle a_1, \ldots, a_n\rangle)$ below.

\begin{defi}\label{def1}
Let $S$ be a formal theory which extends $\Zi$. Fix a G\"odel numbering
of its formulas. Assume there exist recursive functions $\ax$ and $\ded$
such that $\ax$ enumerates the G\"odel numbers of the axioms of $S$
and $\ded$ enumerates all triples $\langle \ulcorner \Ac\urcorner,
\ulcorner \Bc\urcorner, \ulcorner\Cc\urcorner\rangle$ such that $S$
has a deduction rule that infers $\Cc$ from $\Ac$ and $\Bc$ (perhaps
with $\Ac = \Bc$). Also fix recursive functions $f$ and $g$ such that
if $n = \ulcorner \Ac(v_i)\urcorner$ then $f(n,i,k) =
\ulcorner \Ac(\overline{k})\urcorner$ (i.e., all free occurences of $v_i$
are replaced by $\overline{k}$) and $g(n,i) =
\ulcorner (\forall v_i) \Ac(v_i)\urcorner$, where $v_i$ is the $i$th
number variable symbol. If $n$ is not the G\"odel number of a formula,
assume $f(n,i,k) = g(n,i) = 0$.

We define the {\it Tarskian truth theory of $S$}, $\tar(S)$, to be the
theory whose language is the language of $S$ together with one additional
unary relation symbol $T$ and whose non-logical axioms are those of $S$, with
the induction schema extended to the language of $\tar(S)$, together with
the three axioms
\begin{eqnarray*}
&T(\ax(n))&\\
&T(\pi_1(\ded(n))) \wedge T(\pi_2(\ded(n))) \rightarrow T(\pi_3(\ded(n)))&\\
&(\forall k)T(f(n,i,k)) \leftrightarrow T(g(n,i))&
\end{eqnarray*}
and the axiom schema
$$\Ac(v_1, \ldots, v_j) \leftrightarrow
T(\ulcorner \Ac(\overline{v_1}, \ldots, \overline{v_j})\urcorner)$$
for all formulas $\Ac$ in the language of $S$ with no free set variables
and with all free number variables among $v_1, \ldots, v_j$.
\end{defi}

Less rigorously (but perhaps more readably), the three extra axioms
of $\tar(S)$ assert $T(\ulcorner \Ac\urcorner)$ for every axiom $\Ac$
of $S$; $T(\ulcorner \Ac\urcorner) \wedge T(\ulcorner \Bc\urcorner)
\rightarrow T(\ulcorner \Cc\urcorner)$ whenever there is a deduction rule
of $S$ that infers $\Cc$ from $\Ac$ and $\Bc$; and the $\omega$-rule
$(\forall k)T(\ulcorner \Ac(\overline{k})\urcorner) \leftrightarrow
T(\ulcorner (\forall v)\Ac(v)\urcorner)$. The extra axiom schema
is of course just Tarski's truth condition.

The above definition is imprecise in that one really needs to specify not
merely the functions $\ax$, $\ded$, $f$, and $g$ but their codes (and
similarly for the unnamed function used in the final axiom schema to
substitute numerals for variables in G\"odel numbers, though this could
be defined in terms of $f$). Moreover,
these codes must be chosen in a natural way in order to allow the
formalization of proofs that are given below only informally. The most
straightforward way to handle this rigorously would be to make the
definitions explicit as primitive recursive functions for the systems with
which we will be concerned; however, this would be somewhat tedious and
I leave the reader to convince himself that it is possible. This comment
will also apply to similar definitions that appear later.

Note that it is easy to define a provability predicate ${\rm Prov}$ in
terms of $\ax$ and $\ded$ and to show that $\tar(S)$ proves every
instance of the schema $(\forall n)[{\rm Prov}(\ulcorner \Ac(\overline{n})\urcorner)
\rightarrow \Ac(n)]$.

The general concept of truth may or may not be philosophically
problematic, but here we are only using a limited non self-referential
form which I do not think should be controversial, even for use by a
predicativist. Indeed, one could argue that $\tar(S)$ minus its
induction and $\omega$-rule schemas merely formalizes the assertion
that one accepts $S$, and that these schemas are clearly predicatively
legitimate.${}^{17}$

\subsection{Iterated truth theories}\label{sect2g}

We now describe a way to construct iterated families of truth theories.

\begin{defi}\label{def2}
Let $S$ be a theory which extends $\Zi$ and let $\prec$
be a recursive total order on $\omega$. We define the {\it iterated Tarskian
truth theory of $S$ along $\prec$}, $\tar_\prec(S)$, as follows.
Its language is the language of $S$ together with additional unary relation
symbols $\Acc$ and $T_a$ (for each $a \in \omega$).
Its non-logical axioms are the axioms of $S$, with the induction schema
extended to the language of $\tar_\prec(S)$, together with the axiom
$${\rm Prog}(\Acc)$$
(stating progressivity of $\Acc$ with respect to $\prec$). It also has
an additional set of deduction rules whose statement requires some
preparation.

Say that a formula is {\it readable by $T_a$} if it is a formula
of the language of $S$ enriched by the unary relation symbols $T_b$ for
$b \prec a$. Fix a G\"odel numbering of the formulas of $\tar_\prec(S)$
such that both the function $a \mapsto \ulcorner T_a \urcorner$
and the relation $\Rd(a,n)$ indicating that $n$ is the G\"odel
number of a formula readable by $T_a$ are recursive. We also
assume there is a recursive function $\ax$ such that $\ax(a,\cdot)$
enumerates the G\"odel numbers of the axioms of $S$ together with all
logical axioms and the induction schema extended to all formulas readable
by $T_a$, and a recursive function $\ded$ such that $\ded(a,\cdot)$
enumerates the logical deduction rules extended to all formulas
readable by $T_a$ (via triples, in the same manner as in Definition
\ref{def1}). We also require recursive functions $f$ and $g$
satisfying similar modifications of the corresponding conditions in
Definition \ref{def1}; a recursive function $h$ such that
$$h(a,\ulcorner \Ac\urcorner) =
\ulcorner\Ac(v_1, \ldots, v_j)\leftrightarrow
T_a(\ulcorner \Ac(\overline{v_1}, \ldots, \overline{v_j})\urcorner)\urcorner$$
for every formula $\Ac$ in the language of $\tar_\prec(S)$ with no free set
variables and with all free number variables among $v_1, \ldots, v_j$; and
a recursive relation $\Bd$ of $S$ such that $\Bd(\ulcorner\Ac\urcorner)$ holds
if and only if $\Ac$ has no free set variables, for every formula $\Ac$ of
$\tar_\prec(S)$. The extra deduction rules of $\tar_\prec(S)$ then state,
for each $a \in \omega$, that one can infer from $\Acc(\overline{a})$ the
assertions
\begin{eqnarray*}
&T_a(\ax(\overline{a},n))&\\
&T_a(\pi_1(\ded(\overline{a},n))) \wedge T_a(\pi_2(\ded(\overline{a},n)))
\rightarrow T_a(\pi_3(\ded(\overline{a},n)))&\\
&\Rd(\overline{a},n) \rightarrow
\left[(\forall k)T_a(f(\overline{a},n,i,k))
\leftrightarrow T_a(g(\overline{a},n,i))\right]&\\
&\left[b \prec \overline{a} \wedge \Rd(b,n) \wedge \Bd(n)\right]
\rightarrow T_a(h(b,n))&
\end{eqnarray*}
and the assertions
$$\Ac(v_1, \ldots, v_j) \leftrightarrow
T_a(\ulcorner \Ac(\overline{v_1}, \ldots, \overline{v_j})\urcorner)$$
for all formulas $\Ac$ readable by $T_a$ with no free set variables
and with all free number variables among $v_1, \ldots, v_j$. This
completes the definition of $\tar_\prec(S)$.
\end{defi}

The statement $\Acc(a)$ is supposed to signify that one accepts
the truth predicate $T_a$. Thus, this set-up allows the presence
of truth predicates which are not initially known to be intelligible
and can be reasoned with only after some criterion that convinces us
of their legitimacy is satisfied. Specifically, the criterion for
accepting $T_a$ is that we should be assured of the intelligibility
of all formulas to which it applies, which just means that we
should have accepted all prior truth predicates. Thus, the one extra
axiom of $\tar_\prec(S)$ states that $\Acc$ is progressive. The extra
deduction rules allow a predicativist, once he has accepted the truth
predicate $T_a$, to invoke all of the axioms appropriate to that predicate.
So if a predicativist accepts $S$ he should also accept $\tar_\prec(S)$,
for any recursive total order $\prec$ on $\omega$.

The system $\tar_\prec(S)$ is open to the objection that to a limited
extent it allows one to reason with predicates $T_a$ that are not, or not
yet, known to be intelligible. A more stringent and perhaps preferrable
formulation of the theory would ban any use of formulas involving $T_a$
until $\Acc(\overline{a})$ has been proven; this
could be accomplished by deleting all axioms and
all logical deduction rules which involve any truth predicates, and adding
new deduction rules which allow their implementation after an appropriate
acceptability statement has been proven. Similarly, if there are concerns
about interpreting $T_a(n)$ when $n$ is not the G\"odel number of a formula
readable by $T_a$, it is possible to set up a system of distinct G\"odel
numberings, one for each $a \in \omega$, such that every $n \in \omega$ is
the G\"odel number of a formula readable by $T_a$ in the appropriate numbering,
for each $a$. This approach would involve heavy use of a recursive translation
function which relates distinct G\"odel numberings. In any case the net result
would be a slightly more complicated theory with precisely the same
deductive power. In particular, the results presented below would still hold.

\subsection{${\bf Tarski}_{\Gamma_0}^\omega({\bf Z}_1^i)$ and
$\Gamma_0$}\label{sect2h}

I will now present a predicatively valid system
$\tar_{\Gamma_0}^\omega(\Zi)$
which proves well-ordering in a strong sense for notations for every
ordinal less than \Gm{}. The interest of this system is that it shows
that every ordinal less than \Gm{} is predicatively provable. This claim
has been made many times before, but as I pointed out in Section
\ref{sect1}, all previous efforts have crucially involved
impredicative reasoning. Moreover, it will be evident that the
present construction can easily be pushed further to obtain predicative
well-ordering proofs for notations for even larger ordinals. I will describe
stronger systems that go significantly further in \S \ref{sect2i} and
\S \ref{sect2j}.

\begin{defi}\label{def3}
Let $\prec$ be a standard recursive ordering of $\omega$ of order type \Gm{}
whose least elements are $0$ and $1$. We write $\tar_{\Gamma_0}$ for
$\tar_\prec$. Define $\tar_{\Gamma_0}^0(\Zi) = \Zi$ and inductively set
$\tar_{\Gamma_0}^{n+1}(\Zi) = \tar_{\Gamma_0}(\tar_{\Gamma_0}^n(\Zi))$.
Observe that $\tar_{\Gamma_0}^{n+1}(\Zi)$ extends $\tar_{\Gamma_0}^n(\Zi)$.
Let $\tar_{\Gamma_0}^\omega(\Zi)$ be the union of the theories
$\tar_{\Gamma_0}^n(\Zi)$.
\end{defi}

I argued in \S \ref{sect2g} that if a predicativist accepts a theory $S$
then he should also accept $\tar_\prec(S)$ for any recursive total order
$\prec$. In particular, if he accepts $S$ he should also accept
$\tar_{\Gamma_0}(S)$. Granting that he accepts $\Zi$, by iteration
he should accept each $\tar_{\Gamma_0}^n(\Zi)$, and recognizing this
he should also accept $\tar_{\Gamma_0}^\omega(\Zi)$. Now there
is no reason to stop at $\omega$, and by going further we can obtain
predicative well-ordering proofs of larger ordinals. In particular,
$\tar(\tar_{\Gamma_0}^\omega(\Zi))$ proves that a notation for \Gm{}
is well-ordered, yielding the falsification of the \Gm{} thesis promised
in the title and at the beginning of Section \ref{sect1}.

The well-ordering proof is based on the following lemma. Let $a_n$ be the
notation for $\gamma_n$ according to $\prec$. If $a$ is the notation for
$\alpha$ and $b$ is the notation for $\beta$ then let $a \,\dot{+}\, b$
be the notation for $\alpha + \beta$, $\dot{\omega}^a$ the notation for
$\omega^\alpha$, etc. (Let $a\cdot b$ be the notation for $\alpha\beta$.)
We may assume that these are all recursive functions of $a$ and $b$.

\begin{lemma}\label{lemma1}
Let $S$ be a formal theory that extends $\Zi$ and satisfies all assumptions
needed in Definition \ref{def2}.
Then for any $n \in \omega$, $\tar_{\Gamma_0}(S)$ plus the transfinite
induction schema ${\rm TI}(\Ac,a_n)$ for every formula $\Ac$ in its language
proves ${\rm TI}(\Ac,a_{n+1})$ for every formula $\Ac$ in the language of $S$.
\end{lemma}

\begin{proof}
For $n = 0$ one simply carries out the proof of Lemma 1 on page 180 of
\cite{Sch3} within $\tar(S)$. Thus, fix $n \geq 1$. Let $\Bc_{a,b}(m)$ be
the formula
$$\Bc_{a,b}(m) \equiv
(\forall z)\left[\Rd(\dot{\omega}^a\cdot b, z)
\wedge \Bd(z) \rightarrow
T_{\dot{\omega}^a\cdot b}(\ulcorner{\rm Prog} [[z]] \rightarrow
\Jc([[z]], \dot{\phi}_{\overline{a}}(\overline{m}))\urcorner)\right]$$
where
$$\Jc(\Ac,a) \equiv
(\forall y) \left[(\forall x \prec y)\Ac(x)
\rightarrow (\forall x \prec y \,\dot{+}\, a)\Ac(x)\right]$$
and the $[[z]]$ notation indicates that the formula with G\"odel number
$z$ is to be inserted at that point. Note that $\Bc_{a,b}(m)$ is not a
single formula in the language of $\tar_{\Gamma_0}(S)$ because of the
presence of the varying unary relation symbols $T_{\dot{\omega}^a\cdot b}$, but
$$\Cc(a) \equiv
(\forall b)\left[0 \prec b \prec a_n \rightarrow
T_{a_n}(\ulcorner{\rm Prog}_m\,
\Bc_{a,b}(m)\urcorner)\right]$$
is a single formula with parameter $a$.

Since $\Acc$ is progressive, the hypothesis about transfinite induction
yields $\Acc(a_n)$, so that the axioms for $T_{a_n}$ are available.
We will use them to prove in $\tar_{\Gamma_0}(S)$ that $\Cc(a)$ is
progressive over $a \prec a_n$. First, $\Cc(0)$ can be proven by carrying
out, within $T_{a_n}$, the proof of Lemma 1 on page 180 of \cite{Sch3}.
$\Cc(a)$ holds at limit values of $a$ if it holds at all smaller values
by a straightforward verification using the facts that $\phi_a(0) =
\sup_{\tilde{a} \prec a} \phi_{\tilde{a}}(0)$, $\phi_a(m) =
\phi_{\tilde{a}}(\phi_a(m))$ for all $\tilde{a} \prec a$, and $\phi_a(m+1)
= \sup_{\tilde{a} \prec a} \phi_{\tilde{a}}(\phi_a(m)+1)$, all of which
are provable in $\Zi$. At successor stages, we assume $\Cc(a)$ and prove
$\Cc(a\dot{+}1)$ as follows. Fix $0 \prec b \prec a_n$ and $m \in \omega$
and, working within $T_{a_n}$, suppose $\Bc_{a\dot{+}1,b}(m)$ holds; we must
prove $\Bc_{a\dot{+}1,b}(m\dot{+}1)$. This will verify progressivity at
successor stages. $\Bc_{a\dot{+}1, b}(m)$ is proven similarly for $m = 0$,
and it is trivial at limit stages assuming it holds at all previous stages.

The following argument is carried out within $T_{a_n}$. To prove
$\Bc_{a\dot{+}1,b}(m\dot{+}1)$, observe that (by $\Cc(a)$)
$\Bc_{a, r\cdot b}(s)$ is progressive in $s$ for all $r \prec \dot{\omega}$.
Since $\Bc_{a, r\cdot b}(s)$ is readable by
$T_{\dot{\omega}^{a\dot{+}1}\cdot b}$ and has no free set variables,
the hypothesis that $\Bc_{a\dot{+}1,b}(m)$ holds then implies
$$\Jc(\Bc_{a,r\cdot b}(s), \dot{\phi}_{a\dot{+}1}(m))$$
for all $r$. Fixing $j,r \prec \dot{\omega}$ we successively infer
$$\Jc(\Bc_{a,(r\dot{+}j)\cdot b}(s), \dot{\phi}_{a\dot{+}1}(m))$$
and hence
$$\Bc_{a,(r\dot{+}j)\cdot b}(\dot{\phi}_{a\dot{+}1}(m)\dot{+}1);$$
then if $j \succeq 1$, since $\Bc_{a,(r\dot{+}j\dot{-}1)\cdot b}(s)$ is
readable by $T_{\dot{\omega}^a\cdot(r\dot{+}j)\cdot b}$,
$$\Jc(\Bc_{a,(r\dot{+}j\dot{-}1)\cdot b}(s),
\dot{\phi}_a(\dot{\phi}_{a\dot{+}1}(m)\dot{+}1))$$
and hence
$$\Bc_{a,(r\dot{+}j\dot{-}1)\cdot b}(\dot{\phi}_a(\dot{\phi}_{a\dot{+}1}(m)\dot{+}1));$$
and so on, down to
$$\Bc_{a,r\cdot b}(\dot{\phi}_a^j(\dot{\phi}_{a\dot{+}1}(m)\dot{+}1)).$$
Since $j \prec \dot{\omega}$ is arbitrary, this implies
$$\Bc_{a,r\cdot b}(\dot{\phi}_{a\dot{+}1}(m\dot{+}1)),$$
and since this is true for all $r \prec \omega$ (and every formula
readable by $T_{\dot{\omega}^{a\dot{+}1}\cdot b}$ is readable by
$T_{\dot{\omega}^a\cdot r\cdot b}$ for some $r$), we infer
$\Bc_{a \dot{+} 1, b}(m\dot{+} 1)$, as desired.

We conclude that $\Cc(a)$ is progressive, so our hypothesis about
transfinite induction yields $\Cc(a)$ for all $a \prec a_n$. Taking
$b = 1$ and $m = 0$, we infer
$$(\forall a \prec a_n)T_{a_n}(\ulcorner \Bc_{a,1}(0)\urcorner).$$
In particular, for every formula $\Ac$ in the language of $S$ with no
free set variables and one free number variable we have
$$(\forall a \prec a_n)T_{a_n}(\ulcorner{\rm Prog}\, \Ac \rightarrow
\Jc(\Ac,\dot{\phi}_{\overline{a}}(0))\urcorner)$$
and therefore
$${\rm Prog}\, \Ac \rightarrow
(\forall a \prec a_n) \Jc(\Ac,\dot{\phi}_{\overline{a}}(0))$$
and finally (since $a_{n+1} = \phi_{a_n}(0) = \sup_{a \prec a_n} \phi_a(0)$)
$${\rm Prog}\, \Ac \rightarrow (\forall a \prec a_{n+1}) \Ac(a).$$
This proves ${\rm TI}(\Ac,a_{n+1})$ for every formula $\Ac$ in the language
of $S$ with no free set variables and one free number variable. We can reduce
to this case by replacing an arbitrary formula $\Ac$ with the formula
${\rm Prog}(\Ac) \rightarrow \Ac$ and then universally quantifying all
parameters. This yields a formula which is automatically progressive
and we conclude by the above that it holds for all $a \prec a_{n+1}$;
interchanging the order of universal quantifiers then yields
${\rm Prog}(\Ac) \rightarrow (\forall a \prec a_{n+1}) \Ac(a)$.
\end{proof}

This lemma could also be proven by adapting the proof of Theorem 3
in \cite{FefSt}. The point is that one can model the jump hierarchy
used there by using the recursion theorem to find $r \in \omega$ such
that
$$\{r\}(0,m) = \ulcorner\Ac(\overline{m})\urcorner$$
and for all $a \succ 0$
$$\{r\}(a,m) =
\ulcorner (\forall z)[z \prec \overline{a} \rightarrow
\Jc(T_a(\{\overline{r}\}(z,x)),
\dot{\phi}(e(\overline{a}),\overline{m}))]\urcorner,$$
both provably in $\Zi$. Setting $\Bc(a,m) \equiv
T_{a_n\dot{+}1}(\{\overline{r}\}(a,m))$, we can then prove in
$\tar_{\Gamma_0}(S)$ that for all $a \preceq a_n$
$$\Bc(a,m) \leftrightarrow
(\forall z)[z \prec a \rightarrow \Jc(\Bc(z,x), \dot{\phi}(e(a),m))].$$
The argument used to prove Theorem 3 of \cite{FefSt} can then be used to
complete the proof of the lemma. (It can even be slightly simplified
because $\omega^{\gamma_n} = \gamma_n$ for $n \geq 1$ and the condition
$h(a) \preceq z$ is not really necessary, there or here.) However, the proof
given above generalizes better to the situation in \S \ref{sect2i}
and \S \ref{sect2j}.

\begin{theo}\label{theorem1}
$\tar^\omega_{\Gamma_0}(\Zi)$ proves transfinite induction up to any
ordinal less than $\Gamma_0$ for all formulas in its language.
\end{theo}

\begin{proof}
Inductive application of the lemma shows that $\tar_{\Gamma_0}^n(S)$ proves
transfinite induction up to $\gamma_n$ for all formulas of $S$. Applying
this result with $S = \tar_{\Gamma_0}^k(\Zi)$ for arbitrary $n$ and $k$,
and observing that every formula of $\tar^\omega_{\Gamma_0}(\Zi)$ is
a formula of $\tar_{\Gamma_0}^k(\Zi)$ for some $k$, we obtain the desired
result.
\end{proof}

The main point of this proof is that we finesse the induction versus
recursion issue encountered in Section \ref{sect1} by proving a stronger
result, namely that transfinite induction holds for arbitrary formulas
rather than just for sets.

I believe this is the first predicatively valid demonstration of the
well-foundedness of notations for all ordinals less than \Gm{}. Now its
predicative validity requires the validity of three stages of abstraction:
reasoning using Tarskian truth predicates (going from $S$ to $\tar(S)$),
reasoning about which of a sequence of truth predicates are acceptable
(going from $S$ to $\tar_\prec(S)$), and iterating the preceding step
(going from $S$ to $\tar_\prec^\omega(S)$).

If we grant that a predicativist can always pass from $S$ to $\tar_\prec(S)$,
then he should be able to consider the sequence of theories $\tar_\prec^n(S)$
and reason that the validity of each one implies the next, hence they are
all valid and therefore so is $\tar_\prec^\omega(S)$. Thus, if we accept the
second stage then we should also accept the third stage. Conceivably one could
try to make a case that a predicativist can pass to $\tar_\prec^{n+1}(S)$ once
he has actually accepted $\tar_\prec^n(S)$ but he cannot recognize that
this passage is valid in general. This would enable him to accept each
$\tar_\prec^n(S)$ but not $\tar_\prec^\omega(S)$. Presumably the idea
would be that predicativists are capable of making individual judgements
about the acceptability of particular theories but cannot reason about the
acceptability of theories in general. However, it would seem that one could
with equal justice replace ``theories'' with ``truth predicates'' in this
assertion, yielding the claim that predicativists can make individual
judgements about whether to accept particular truth predicates but cannot
reason about such questions abstractly, which would render $\tar_\prec(S)$
illegitimate. In other words, we have just as much reason to go from the
second stage to the third as we do to go from the first stage to the second.

On the other hand, if we reject passage from the first stage (accepting
$\tar(S)$) to the second (accepting $\tar_\prec(S)$) I think it would have
to be on the grounds just suggested, i.e., that predicativists are capable
of making individual judgements about the acceptability of particular truth
predicates but cannot reason abstractly about their acceptability in
general. But then we could just as well replace ``truth predicates'' with
``statements'' and argue that predicativists can make individual judgements
about whether to accept particular statements but cannot reason about such
questions abstractly. This would render $\tar(S)$ illegitimate. But not
only that; as I pointed out in \S \ref{sect1d} (a), this line of argument,
if accepted, would actually prevent a predicativist from recognizing modus
ponens, or really any deduction rule, as a general principle. In fact it
ought to forbid any use of statements involving variables of any kind
since these already imply an ability to reason hypothetically about the
truth of such a statement on all possible substitutions of values for
the variables, which evidently requires some abstract sense of truth
or acceptability. Ultimately we would be left only with the ability to
make concrete numerical assertions. The point is that
{\it there is a smooth progression in reasoning at successively
higher levels of abstraction that takes one from Peano arithmetic up to
$\tar_{\Gamma_0}^\omega(\Zi)$ and beyond}, so that any attempt to cut
this progression off at some point is bound to appear arbitrary. At any
rate there is no better reason to cut it off at \Gm{} than anywhere else.

The element of truth to the objection is that the general concept of
statements ``being predicatively acceptable'' is impredicative when there
is no limitation on the domain of discussion, because it would be circular
to talk about the predicative acceptability of statements which themselves
involve the concept of predicative acceptability. Similarly, the ideas
of truth predicates being predicatively acceptable or of theories being
predicatively acceptable are impredicative as unrestricted general concepts.
But these are just versions of the fact noted in \S \ref{sect2e} that
the general concept of ``being intelligible'' is not itself intelligible.
In each of these cases, if attention is restricted to a well-defined
previously grasped domain, I see no impredicativity.

\subsection{${\bf Tarski}_{\kappa^\kappa}^\omega({\bf Z}_1^i)$ and
$\phi_{\Omega^2}({\bf 0})$}\label{sect2i}

In the last section we considered a formal system, $\tar_{\Gamma_0}(\Zi)$,
in which we were able to reason about the (predicative) acceptability of
a hierarchy of truth predicates by means of an additional
predicate $\Acc$. The $\tar_{\Gamma_0}$ construction was then iterated
$\omega$ times. By systematizing the process of iterating constructions
involving acceptability predicates we can access ordinals well beyond
\Gm{}. I will next illustrate this claim by describing a predicatively
valid formal system that proves well-ordering statements for notations
for all ordinals less than the Ackermann ordinal $\phi_{\Omega^2}(0)$.
In \S \ref{sect2j} I will sketch a way to carry the construction further
and get a predicative well-ordering proof for the ``small'' Veblen ordinal
$\phi_{\Omega^\omega}(0)$.

The development is similar to that in \S \ref{sect2h} and will be presented
here in slightly less detail. Let $\kappa = \phi_{\Omega^\Omega}(0)$ and fix
a notation system for $\kappa^\kappa$ (e.g., see the introduction to
\cite{Mill}). In the following I will identify ordinals with their
notations and I will use $\alpha$, $\beta$, $\gamma$ to range over
ordinals $\prec \kappa$ and $a$, $b$, $c$ to range over ordinals
$\prec \kappa^\kappa$. Every nonzero $a$ can be uniquely written in the
form $a = \kappa^{\alpha_1}\beta_1 + \cdots + \kappa^{\alpha_n}\beta_n$
such that $\alpha_1 \succeq \cdots \succeq \alpha_n$; $\alpha_i = \alpha_{i+1}$
implies $\beta_i \succeq \beta_{i+1}$; and each $\beta_i$ is either 1 or
a limit ordinal. Let $h(a) = \kappa^{\alpha_1}\beta_1 + \cdots +
\kappa^{\alpha_{n-1}}\beta_{n-1}$. We define the {\it canonical sequence
associated to $a$} to be $\{h(a) + \kappa^{\alpha_n}\gamma: \gamma \prec
\beta_n\}$ if $\beta_n$ is a limit ordinal; if $\beta_n = 1$ and $\alpha_n
= 0$ then it is the single element $\{h(a)\}$; if $\beta_n = 1$ and
$\alpha_n$ is a limit ordinal it is
$\{h(a) + \kappa^\gamma: \gamma \prec \alpha_n\}$; and if $\beta_n = 1$
and $\alpha_n = \tilde{\alpha}_n + 1$ then it is $\{h(a) +
\kappa^{\tilde{\alpha}_n}\gamma: \gamma \prec \kappa\}$. In the last
case ($\beta_n = 1$ and $\alpha_n$ a successor) we say that $a$ is of
{\it type 1}, and otherwise it is of {\it type 0}. (Cf.\ Definition 1
of \cite{Ger}.) We consider $0$ to be of type 0 and we let its canonical
sequence be empty. Let $\Typ_0$ be a formula such that $\Typ_0(a)$ holds
if and only if $a$ is of type 0, and let $\Seq$ be a formula such that
$\Seq(x,a)$ holds if and only if $x$ belongs to the canonical sequence
associated to $a$.

\begin{defi}\label{def4}
Let $S$ be a theory which extends $\Zi$. We define $\tar_{\kappa^\kappa}(S)$
as follows. It language is the language of $S$ together with a unary relation
symbol $\Acc$ and two families of unary relation symbols $T_a$ and $\Acc_a$.
In this setting a formula is {\it readable by $T_a$} if it is a formula of
the language of $S$ enriched by the unary relation symbols $T_b$ and $\Acc_b$
for all $b \prec a$. The non-logical axioms are the axioms of $S$, with the
induction schema extended to the larger language, together with an axiom
which states that for any $a$, if $a$ is of type 0 with associated canonical
sequence $(x_\gamma)$ then $(\forall\gamma)\Acc(x_\gamma) \rightarrow \Acc(a)$,
and if $a$ is of type 1 with canonical sequence $(x_\gamma)$ then
${\rm Prog}_\gamma\, \Acc(x_\gamma) \rightarrow \Acc(a)$. We also have, for
each $a$, a family of deduction rules allowing inference from the premise
$\Acc({\overline{a}})$ of the following statements:
\medskip

\noindent (I) Axioms for $T_a$: the same as in Definition \ref{def2}, i.e.,
\begin{eqnarray*}
&T_a(\ax(\overline{a},n))&\\
&T_a(\pi_1(\ded(\overline{a},n))) \wedge T_a(\pi_2(\ded(\overline{a},n)))
\rightarrow T_a(\pi_3(\ded(\overline{a},n)))&\\
&\Rd(\overline{a},n) \rightarrow
\left[(\forall k)T_a(f(\overline{a},n,i,k))
\leftrightarrow T_a(g(\overline{a},n,i))\right]&\\
&\left[b \prec \overline{a} \wedge \Acc_a(b) \wedge \Rd(b,n)
\wedge \Bd(n)\right] \rightarrow T_a(h(b,n))&\\
&\Ac(v_1, \ldots, v_j) \leftrightarrow
T_a(\ulcorner \Ac(\overline{v_1}, \ldots, \overline{v_j})\urcorner),&
\end{eqnarray*}
with the premise $\Acc_a(b)$ added to the fourth axiom. The functions
appearing in these axioms are defined analogously to those in Definition
\ref{def2}, and the same condition is placed on $\Ac$ in the final schema.
\medskip

\noindent (II) Axioms for $\Acc_a$:
\begin{eqnarray*}
&\Seq(b,\overline{a}) \rightarrow \Acc_a(b)&
\!\!\!\!\!\!\!\!\!\!\!\!\!\!\!\hbox{(if $a$ is of type 0)}\\
&{\rm Prog}_b\left[\Seq(b,\overline{a}) \rightarrow \Acc_a(b)\right]&
\!\!\!\!\!\!\!\!\!\!\!\!\!\!\!\hbox{(if $a$ is of type 1)}\\
&[b \prec \overline{a} \wedge \Acc_a(b) \wedge \Typ_0(b)] \rightarrow
T_a\left(\ulcorner(\forall c)[\Seq(c,b) \rightarrow
\Acc_b(c)]\urcorner\right)&\\
&[b \prec \overline{a} \wedge \Acc_a(b) \wedge \neg\Typ_0(b)] \rightarrow
T_a\left(\ulcorner{\rm Prog}_c [\Seq(c,b) \rightarrow
\Acc_b(c)]\urcorner\right)&\\
&\left[c \prec b \prec \overline{a} \wedge \Acc_a(b)\right] \rightarrow
\left[\Acc_a(c) \leftrightarrow T_a(\ulcorner\Acc_b(c)\urcorner)\right]&\\
&b \prec \overline{a} \rightarrow \left[\Acc(b) \leftrightarrow
\Acc_a(b)\right]&
\end{eqnarray*}
This completes the definition of $\tar_{\kappa^\kappa}(S)$.

As in Definition \ref{def3} we now inductively define
$\tar_{\kappa^\kappa}^0(\Zi) = \Zi$ and $\tar_{\kappa^\kappa}^{n+1}(\Zi) =
\tar_{\kappa^\kappa}(\tar_{\kappa^\kappa}^n(\Zi))$, and we
let $\tar_{\kappa^\kappa}^\omega(\Zi)$ be the union of the
theories $\tar_{\kappa^\kappa}^n(\Zi)$.
\end{defi}

The system $\tar_{\kappa^\kappa}^\omega(\Zi)$ is predicatively justified
by taking one step up in abstraction beyond $\tar_{\Gamma_0}(\Zi)$.
There we had a theory involving a hierarchy of truth predicates and we
formally reasoned about their acceptability. Here we have a hierarchy
of acceptability predicates, each of which allows us to reason about
the acceptability of truth and acceptability predicates of lower degree,
and about whose acceptability we are able to formally reason. Intuitively,
if $a$ is of type 0 then the acceptability predicate at level $a$ is
supposed to affirm the acceptability of all levels belonging to the
canonical sequence associated to $a$, and if $a$ is of type 1 then it
is supposed to affirm progressivity of the acceptability of the levels
belonging to the canonical sequence associated to $a$.

Just as with Definition \ref{def2} it is possible to formulate a
stricter definition which would disallow any use of $T_a$ and $\Acc_a$
until after $\Acc(\overline{a})$ has been proven. We could also set
up a family of distinct G\"odel numberings so that truth predicates
could only refer to formulas readable by them, and similarly we could
set up a family of distinct orderings of $\omega$, one of order type
$a$ for each $a$, so that acceptability predicates could only refer to
prior truth and acceptability predicates. As before, these changes
would be merely cosmetic and would not affect the strength of the theory.

The well-ordering proof is based on two lemmas.

\begin{lemma}\label{lemma2}
Let $S$ be a formal theory that extends $\Zi$ and satisfies all assumptions
needed in Definition \ref{def4}. Then $\tar_{\kappa^\kappa}(S)$ proves that
the statement
$$\Ac(\alpha) \equiv
(\forall a)\left[\Acc(a) \rightarrow \Acc(a + \kappa^\alpha)\right]$$
is progressive in $\alpha$.
\end{lemma}

\begin{proof}
For any $a$ the canonical sequence associated to $a + 1$ is $\{a\}$, so an
axiom of $\tar_{\kappa^\kappa}(S)$ asserts that $\Acc(a) \rightarrow
\Acc(a+1)$. This shows that $\Ac(0)$ is provable in $\tar_{\kappa^\kappa}(S)$.
At successor stages, suppose $\Ac(\alpha)$
holds and, reasoning in $\tar_{\kappa^\kappa}(S)$, deduce that for
any $a$ satisfying $\Acc(a)$ the statement $\Acc(a + \kappa^\alpha\beta)$ is
progressive in $\beta$. This yields $\Acc(a) \rightarrow
\Acc(a + \kappa^{\alpha + 1})$, and we infer $\Ac(\alpha + 1)$.
Finally, suppose $\alpha$ is a limit and we have $\Ac(\gamma)$ for
all $\gamma \prec \alpha$. Then for any $a$ such that $\Acc(a)$ holds
we have $\Acc(a + \kappa^\gamma)$ for all $\gamma \prec \alpha$, and
this implies $\Acc(a + \kappa^\alpha)$. So we infer $\Ac(\alpha)$.
\end{proof}

Let $\delta_0 = 1$ and $\delta_{n+1} = \phi_{\Omega\cdot \delta_n}(0)$,
so that $\phi_{\Omega^2}(0) = \sup_{n \in \omega} \delta_n$. Note that
$\omega\kappa = \kappa$, so $(\kappa\omega)^\alpha$ equals $\kappa^\alpha$
if $\alpha$ is a limit and it equals $\kappa^\alpha\omega$ if $\alpha$ is
a successor.

\begin{lemma}\label{lemma3}
Let $S$ be a formal theory that extends $\Zi$ and satisfies all assumptions
needed in Definition \ref{def4}. Then $\tar_{\kappa^\kappa}(S)$ plus the
transfinite induction schema ${\rm TI}(\Ac,\delta_n)$ for every formula
$\Ac$ in its language proves ${\rm TI}(\Ac,\delta_{n+1})$ for every
formula $\Ac$ in the language of $S$.
\end{lemma}

\begin{proof}
The technique is similar to that used in the proof of Lemma \ref{lemma1}.
Define
$$\Bc_{\alpha,b}(\mu) \equiv
(\forall z)\left[\Rd((\kappa\omega)^\alpha b, z)
\wedge \Bd(z) \rightarrow
T_{(\kappa\omega)^\alpha b}(\ulcorner{\rm Prog} [[z]] \rightarrow
\Jc([[z]], \phi_{\Omega\alpha}(\mu))\urcorner)\right]$$
and
$$\Cc(\alpha) \equiv
(\forall b)\left[0 \prec b \prec \kappa^{\delta_n} \wedge \Typ_0(b)
\wedge \Acc_{\kappa^{\delta_n}}((\kappa\omega)^\alpha b) \rightarrow
T_{\kappa^{\delta_n}}(\ulcorner{\rm Prog}_\mu\,
\Bc_{\alpha,b}(\mu)\urcorner)\right].$$

By Lemma \ref{lemma2} and the transfinite induction hypothesis we obtain
$\Acc(a) \rightarrow \Acc(a + \kappa^\alpha)$ for all $\alpha \prec \delta_n$.
In particular, $\Acc(\kappa^\alpha)$ holds for all $\alpha \prec \delta_n$,
and this implies $\Acc(\kappa^{\delta_n})$. We also obtain
$\Acc(a) \rightarrow \Acc(a + \kappa^\alpha r)$ for all $\alpha \prec
\delta_n$ and all $r \prec \omega$, which implies $\Acc(a) \rightarrow
\Acc(a + \kappa^\alpha\omega)$ for all $\alpha \prec \delta_n$.

We claim that $\Cc(\alpha)$ is progressive over $\alpha \prec \delta_n$.
$\Cc(0)$ is again essentially Lemma 1 on page 180 of \cite{Sch3}. Next,
let $\alpha$ be a limit and suppose $\Cc(\beta)$ holds for all $\beta
\prec \alpha$. Fix $0 \prec b \prec \kappa^{\delta_n}$ of type 0 and
suppose $\Acc_{\kappa^{\delta_n}}((\kappa\omega)^\alpha b)$. Then for
every $\beta \prec \alpha$ and every $x$ in the canonical sequence
associated to $b$, letting $y = (\kappa\omega)^{\tilde{\alpha}}x + 1$
where $\beta + \tilde{\alpha} = \alpha$, we have
$\Acc_{\kappa^{\delta_n}}((\kappa\omega)^\beta y)$ and
so $\Cc(\beta)$ implies progressivity in $\mu$ of the assertion that
$T_{(\kappa\omega)^\beta y}(\ulcorner {\rm Prog} [[z]] \rightarrow
\Jc([[z]], \phi_{\Omega\beta}(\mu))\urcorner)$ holds for all appropriate
$z$. Since every formula readable by $T_{(\kappa\omega)^\alpha b}$ is
readable by $T_{(\kappa\omega)^\beta y}$ for sufficiently large $x$
and $\beta$ and $\phi_{\Omega\alpha}$ enumerates the common values of
$\phi_{\Omega\beta}$ over $\beta \prec \alpha$, this implies the
desired conclusion.

Finally, suppose $\Cc(\alpha)$ holds for some $\alpha \prec \delta_n$; we must
verify $\Cc(\alpha+1)$. To do this fix $0 \prec b \prec \kappa^{\delta_n}$
of type 0 and suppose $\Acc_{\kappa^{\delta_n}}((\kappa\omega)^{\alpha+1}b)$.
We may assume $b$ is a successor, $b = \tilde{b} + 1$, since progressivity of
$\Bc_{\alpha + 1, b}(\mu)$ for these values implies progressivity for type 0
limits by an argument similar to the one used in the previous paragraph.
Now, as in the proof of Lemma \ref{lemma1}, progressivity at limit values of
$\mu$ follows from continuity of $\phi_{\Omega(\alpha+1)}$, and progressivity
at $\mu = 0$ is proven in the same way as progressivity at successor values
of $\mu$. Therefore we consider the successor case. Fix $\mu \prec \kappa$
and, working within $T_{\kappa^{\delta_n}}$, suppose
$\Bc_{\alpha + 1, \tilde{b} + 1}(\mu)$ holds; we must prove
$\Bc_{\alpha + 1, \tilde{b} + 1}(\mu + 1)$.

Let $\tilde{\gamma}_0 = \phi_{\Omega(\alpha+1)}(\mu) + 1$ and inductively
define $\tilde{\gamma}_{n+1} = \phi_{\Omega\alpha + \tilde{\gamma}_n}(0)$,
so that
$\phi_{\Omega(\alpha+1)}(\mu + 1) = \sup_{n \prec \omega} \tilde{\gamma}_n$.
In order to verify $\Bc_{\alpha + 1, \tilde{b}+1}(\mu + 1)$, it will therefore
suffice to check $T_{(\kappa\omega)^{\alpha+1}b}(\ulcorner{\rm Prog}[[z]]
\rightarrow \Jc([[z]], \tilde{\gamma}_n)\urcorner)$ for
arbitrary $n$ (and appropriate $z$). Working within $T_{\kappa^{\delta_n}}$,
we can now mimic the proof of Lemma \ref{lemma1} to show that transfinite
induction up to $\tilde{\gamma}_n$ for all formulas readable by
$T_{(\kappa\omega)^{\alpha+1}\tilde{b} + (\kappa\omega)^\alpha\kappa(r+1)}$
implies transfinite induction
up to $\tilde{\gamma}_{n+1}$ for all formulas readable by
$T_{(\kappa\omega)^{\alpha+1}\tilde{b} +(\kappa\omega)^\alpha\kappa r}$. By
the same argument used in the proof of Theorem \ref{theorem1} we then get
transfinite induction up to any $\tilde{\gamma}_n$ for formulas readable by
any $T_{(\kappa\omega)^{\alpha+1}\tilde{b} + (\kappa\omega)^\alpha\kappa r}$,
i.e., for all formulas readable by $T_{(\kappa\omega)^{\alpha+1}b}$. In
this argument $\tilde{\gamma}_n$ replaces $a_n$,  the expression $\omega^a b$
is modified to
$$(\kappa\omega)^{\alpha+1}\tilde{b} + (\kappa\omega)^\alpha\kappa r +
(\kappa\omega)^\alpha \omega^{\tilde{\alpha}}\tilde{\beta}$$
(where $a$ and $b$ are replaced by $\tilde{\alpha}$ and $\tilde{\beta}$),
$\phi_a$ is modified to $\phi_{\Omega\alpha + \tilde{\alpha}}$, and the
hypothesis that $\Cc(\alpha)$ holds replaces use of Lemma 1 on page 180
of \cite{Sch3}; otherwise the argument is identical. We conclude that
$\Cc(\alpha)$ is progressive.

The remainder of the proof is similar to the proof of Lemma \ref{lemma1}.
\end{proof}

\begin{theo}\label{theorem2}
$\tar^\omega_{\kappa^\kappa}(\Zi)$ proves transfinite induction up to any
ordinal less than $\phi_{\Omega^2}(0)$ for all formulas in its language.
\end{theo}

\begin{proof}
By reasoning identical to that used in the proof of Theorem \ref{theorem1}.
\end{proof}

\subsection{${\bf Tarski}_{\lambda^{\lambda^\omega}}^\omega({\bf Z}_1^i)$ and
$\phi_{\Omega^\omega}({\bf 0})$}\label{sect2j}

It is possible to strengthen the construction of \S \ref{sect2i} to
obtain a predicative well-ordering proof for notations for all
ordinals less than the ``small'' Veblen ordinal $\phi_{\Omega^\omega}(0)$.
I will only sketch the construction and I will omit all proofs.

Let $\lambda = \phi_{\Omega^{\Omega^\omega}}(0)$ and fix a notation
system for $\lambda^{\lambda^\omega}$. Every ordinal $a$ less than
$\lambda^{\lambda^\omega}$ can be uniquely written in the form
$a = \lambda^{\alpha_1}\beta_1 + \cdots + \lambda^{\alpha_n}\beta_n$ with
each $\alpha_i \prec \lambda^\omega$; each $\beta_i \prec \lambda$;
$\alpha_1 \succeq \cdots \succeq \alpha_n$; $\alpha_i = \alpha_{i+1}$
implies $\beta_i \succeq \beta_{i+1}$; and each $\beta_i$ either a
limit or $1$. Write $\alpha_n = \lambda^{\gamma_1}\delta_1 + \cdots
+ \lambda^{\gamma_m}\delta_m$ with similar conditions on the $\gamma_i$
and $\delta_i$ (but now each $\gamma_i \prec \omega$). If $\beta_n$ is
a limit or if $\beta_n = 1$ and $\delta_m$ is a limit then we say that
$a$ is of type 0; if $\beta_n = \delta_m = 1$ then $a$ is of type
$\gamma_m + 1$. Canonical associated sequences are defined just as before
in the type 0 case, for type 1 we take the canonical associated
sequence to be
$$\{h(a) + \lambda^{\lambda^{\gamma_1}\delta_1 + \cdots
+ \lambda^{\gamma_{m-1}}\delta_{m-1}} \cdot \alpha: \alpha \prec \lambda\},$$
and for higher types we take it to be
$$\{h(a) + \lambda^{\lambda^{\gamma_1}\delta_1 + \cdots
+ \lambda^{\gamma_{m-1}}\delta_{m-1} + \lambda^{\gamma_m - 1}\alpha}:
\alpha \prec \lambda\},$$
where $h(a) = \lambda^{\alpha_1}\beta_1 + \cdots +
\lambda^{\alpha_{n-1}}\beta_{n-1}$.
The $\tar_{\lambda^{\lambda^\omega}}(S)$ construction goes like the
$\tar_{\kappa^\kappa}(S)$ construction, with the following elaboration.
Define a {\it bounded jumpability} predicate
$$\Jc^1_{a,b}(\Ac,c) \equiv
(\forall a \preceq x \prec b) \left[\Ac(x) \rightarrow \Ac(x+c)\right]$$
and inductively define a {\it bounded $k$-th order jumpability} predicate
$\Jc^k_{a,b}(\Ac,c)$ by
$$\Jc^{k+1}_{a,b}(\Ac,c) \equiv
(\forall a \preceq x \prec b) \left[\Jc^k_{a,b}(\Ac,x)
\rightarrow \Jc^k_{a,b}(\Ac,x + c)\right].$$
Accepting level $a$ of the truth and acceptability hierarchy affiliated to
$\tar_{\lambda^{\lambda^\omega}}(S)$ should then entail accepting level
$h(a)$ as well as: if $a$ is
type 0, accepting all levels in its associated sequence; if $a$ is
type 1, accepting progressivity of the acceptability of the levels
in its associated sequence; and if $a$ is of type $k \geq 2$,
accepting
$${\rm Prog}_{\alpha \prec \lambda}
\Jc^{k-1}_{h(a),a}(\Acc, \lambda^{\lambda^{\gamma_1}\delta_1 + \cdots
+ \lambda^{\gamma_{m-1}}\delta_{m-1} + \lambda^{\gamma_m - 1}\alpha}).$$
The axiomatization of
$\tar_{\lambda^{\lambda^\omega}}(S)$ is similar to the axiomatization
of $\tar_{\kappa^\kappa}(S)$, with natural modifications to accomodate
higher types. We define $\tar_{\lambda^{\lambda^\omega}}^\omega(\Zi)$
by induction as usual. By a more complicated but not essentially different
argument from the proof of Theorem \ref{theorem2}, we obtain the following
result.

\begin{theo}\label{theorem3}
$\tar^\omega_{\lambda^{\lambda^\omega}}(\Zi)$ proves transfinite induction
up to any ordinal less than $\phi_{\Omega^\omega}(0)$ for all formulas in
its language.
\end{theo}

In this proof the analog of Lemma \ref{lemma2} asserts that for each
$k \prec \omega$, $\tar_{\lambda^{\lambda^\omega}}(S)$ proves
$\Acc(\lambda^{\lambda^k})$. Then the analog of Lemma \ref{lemma3}
asserts that $\tar_{\lambda^{\lambda^\omega}}(S)$ plus transfinite
induction up to $\tilde{\delta}_n$ for all formulas in its language
proves transfinite induction up to $\tilde{\delta}_{n+1}$ for all
formulas in the language of $S$, where $\tilde{\delta}_1 = 1$ and
$\tilde{\delta}_{n+1} = \phi_{\Omega^k \delta_n}(0)$. (This proof
is carried out using truth and acceptability at level
$\lambda^{\lambda^{k-1}}\omega$.)

As with $\tar^\omega_{\Gamma_0}(\Zi)$ we can go one step further
and conclude that the Veblen ordinal is itself predicatively provable.
Although $\phi_{\Omega^\omega}(0)$ is not terribly large among the scale
of proof-theoretically important countable ordinals, this result is still
significant. Probably the most celebrated example of an allegedly
impredicative mainstream theorem, Kruskal's theorem (see, e.g., \cite{Gal}),
is now seen to be predicatively justified. It is equivalent over a weak base
system to the well-ordering of a notation for $\phi_{\Omega^\omega}(0)$
\cite{RW}. 

It should not be difficult to strengthen Theorem \ref{theorem3} so as
to prove the well-foundedness of a notation for the ``large'' Veblen
ordinal $\phi_{\Omega^\Omega}(0)$. But I expect that substantially
larger ordinals can be accessed using predicative methods. This raises
the possibility of a version of Hilbert's program in which theories are
justified via predicative, rather than finitary or intuitionistic,
consistency proofs. The preceding results indicate that this program
is interesting, substantial, and open to exploration. Moreover, if
predicativism given the natural numbers --- or, as I prefer to call it,
{\it mathematical conceptualism} \cite{W1} --- is {\it right}, this program
is of fundamental significance for the foundations of mathematics.

\bigskip
\bigskip

\footnotesize{

\noindent 1. Other concerns may include the awkwardness of predicative
systems in practice and a sense that they are philosophically, as
opposed to mathematically, too limiting. I address these issues in
separate papers \cite{W1,W2}.

\noindent 2. Yet another idea is to assert that we merely {\it believe}
that $A$ could come to accept every statement in $\Sc$, but we do not
{\it know} this. If so, it is possible that $A$ could indeed share this
belief, but without sufficient certainty to allow him to go beyond $\Sc$.

Whether this tactic could work depends on exactly why we have reservations
about what $A$ can accept. It is no good, for example, to say that we are
not sure what an ideal predicativist can accept because there is more than
one version of predicativism, since we can hardly assume that he is unable
to decide which version he prefers. Nor would the argument hold up
if our uncertainty were caused by not knowing whether $A$ could accept
some specific principle $P$, as this would be tantamount to $A$
failing to decide between two versions of his theory.

The argument might succeed if there were an infinite sequence
of principles $(P_n)$ each one of which $A$ might accept or reject.
A case could then be made that we cannot demand that he make a
simultaneous decision on the validity of every $P_n$. However, this
now seems to be a version of the idea that $A$ can accept each
$P_n$ individually but not the entire sequence $(P_n)$, which
is the sort of claim I address in the main text.

\noindent 3. The $H$ system in \cite{Fef1} follows this description
precisely. Systems of ramified analysis like $\Sigma$ and $R$ are a
little more complicated in that each $S_a$ has its own set variables
$X^a$, and legal formulas of $S_a$ must contain only set variables
$X^b$ with $b \preceq a$. These systems are formally more complicated
than $H$, but they are supposed to more transparently model the intuition
of a predicative universe which is only available in stages.

\noindent 4. According to the proof sketched in \cite{Fef1} that
$\Gamma_0 \leq {\overline{\rm Aut}}(S)$, we can find $r \in \omega$
which is the G\"odel number of a recursive
function $\{r\}$ with $\{r\}(n)$ a notation for $\gamma_n$,
$\{r\}(n) <_\Oc \{r\}(n+1)$, and $S_0 \vdash (\forall n)\,
{\rm Prov}_{S_{\{\overline{r}\}(n)}}
(\ulcorner I(\{\overline{r}\}({\overline{n+1}}))\urcorner)$.
Letting $\Ac(n) \equiv I(\{\overline{r}\}(n+1))$ and substituting $\{\overline{r}\}(n)$
for $a$ in ($**$), a simple induction argument yields
$(\forall n)\, I(\{\overline{r}\}(n))$ (note that $S_0$ supports 
complete induction),
from which we deduce $I(\overline{a})$ with $a = 3\cdot 5^r$.

(The expression $\Ac(\{\overline{r}\}(x))$ should be understood as an
abbreviation of a formula which asserts that there exists $y$ such
that $\{r\}(x) = y$ and $\Ac(y)$. Alternatively, we can assume a language
that contains symbols for all primitive recursive functions and reword the
arguments --- here and below --- to ensure that all recursive functions in
use are actually primitive recursive.)

\noindent 5. In \cite{Fef8} Feferman refers to ``the argument that
the characterization of predicativity requires one to go beyond
predicative notions and principles'' (\cite{Fef8}, footnote 6),
which sounds like it could be a version of the general objection
of \S \ref{sect1c}. However, his response (``But the
predicativist $\ldots$'', p.\ 316) seems aimed merely at showing that
the set of all predicatively provable ordinals is not a predicatively
valid set, a view that I agree with (though not for the reason given
there). This should not prevent a predicativist from understanding
the assertion that every $a_n$ is an ordinal notation,
in the notation of \S \ref{sect1c}.

One could possibly make an argument that the statement $(\forall n)\, I(a_n)$
cannot even be predicatively recognized as meaningful, let alone true, on
the ground that the general concept of well-ordering is not available to a
predicativist. Perhaps this is the point of the comment in \cite{Fef8}.
Presumably the idea would be that each $I(a_n)$ can only be understood as
a sensible assertion {\it once it is proven} and not before. This seems
like a difficult position to defend, but in any case it would void the main
argument because if one did accept that some theorems of $S_{a_n}$
cannot be recognized as meaningful until they are actually proven, this
would invalidate any use of reflection principles in the first place.

\noindent 6. Note that the final $\Lc$ in rule V, predicate substitution
(\cite{Fef3}, p.\ 78), should be $\Lc_\exists$.

\noindent 7. Of course, the validity of the functional generating
procedure hinges on the validity of $\exists/P$, so it may be significant
that Feferman refers to ``the correctness of $P$'' and not ``the correctness
of $P$ in conjunction with $\exists/P$''. This goes back to the question
raised in \S \ref{sect1f} (a) about whether predicativists can trust theorems
proven in $\exists/P$, and if not, why it makes sense for them to use this
system at all.

\noindent 8. For the argument to work we have to be able to imagine
someone who can think, for example,

\begin{quotation}
\noindent whenever it is the case that for every number $n$ and
any $\vec{\alpha}$ there exists a $\beta$ satisfying
$\Ac(\vec{\alpha}, n, \beta)$, for any $\vec{\alpha}$ these
$\beta$'s can be unified into a single $\gamma$ satisfying
$\Ac(\vec{\alpha}, n, \gamma_n)$ for all $n$
\end{quotation}

\noindent but who cannot think

\begin{quotation}
\noindent whenever for any $\vec{\alpha}$ a unique $\beta$ exists
satisfying $\Ac(\vec{\alpha}, \beta)$, I can introduce a
functional symbol $F$ such that $\Ac(\vec{\alpha}, F(\vec{\alpha}))$
holds for any $\vec{\alpha}$,
\end{quotation}

\noindent yet who {\it can} think

\begin{quotation}
\noindent I can introduce a
functional symbol $F$ such that $\Ac(\vec{\alpha}, F(\vec{\alpha}))$
holds for any $\vec{\alpha}$
\end{quotation}

\noindent once he has actually proven, for any particular $\Ac$, the
existence for any $\vec{\alpha}$ of a unique $\beta$ satisfying
$\Ac(\vec{\alpha}, \beta)$. This combination of abilities and deficits
strikes me as incoherent.

\noindent 9. Another questionable point in the ``too weak'' category
is the restriction on allowed types in (\cite{FefSt}, p.\ 81). I do
not understand the justification given there, and a corresponding
restriction is not made in the system sketched in \cite{Fef6}. In
light of footnote 2 of \cite{FefSt}, this raises the question
whether \UNFA{} as described in \cite{Fef6} really does have
proof-theoretic ordinal \Gm{}.

\noindent 10. On the other hand, the minimality property of $LFP$
(Ax 4 (ii), p.\ 79) is never used in \cite{FefSt}, so this axiom
could be eliminated without affecting the proof-theoretic strength
of \UNFA{}. The existence of not necessarily minimal fixed points
might be predicatively justifiable if intuitionistic logic is used;
see the discussion of inductively defined classes at the beginning
of \S \ref{sect2f}. However, this more careful analysis also reveals
a fundamental impredicativity in using schematic variables, a point
not discussed above; see \S \ref{sect2e}.

\noindent 11. I should point out that Ferferman has in several places
openly called attention to impredicative aspects of various of his systems.
The impredicativity of the autonomous systems
is commented on in (\cite{Fef3}, p.\ 85), (\cite{Fef4}, p.\ 3),
and elsewhere. (``the well-ordering statement $\ldots$ on the face
of it {\it only impredicatively justifies} the transfinite
iteration of accepted principles up to $a$.'' ``$\ldots$ {\it prima facie}
impredicative notions such as those of ordinals or well-orderings.'')
The impredicativity of $P + \exists/P$ is noted in (\cite{Fef3}, p.\ 92).
(``In $P$ we think of `$X$' as ranging over predicates recognized to
have a definite meaning; this would not seem to admit the properties
expressed by formulas of $\Lc_\exists$.'') The impredicativity
of \RefP{} is noted in (\cite{Fef4}, pp.\ 41 and 42). (``one may
question substituting possibly indeterminate formulas $\ldots$
this seems to me to be the weakest point of the case for reflective
closure having fundamental significance.'' ``this may involve some
equivocation between the notions of being definite $\ldots$ and
being determinate'')

\noindent 12. (Perhaps also assuming that we know $Y$ to be countable so
that the subset can be extracted using principle (iii); but see below.)

\noindent 13. I am deliberately avoiding the question of what sets
``really are''. The argument in this paragraph suggests a
quasi-physical conception according to which one could imagine
actually manipulating the elements of a set. I see nothing wrong
with this sort of conception, at least for subsets of $\omega$, but
it is not essential for what follows. The important point is that
sets of numbers, whatever one takes them to be, should in principle
always be unequivocally recognizable as such. This ought to be true
on any reasonable conception (but probably would not be true, for
example, if one identified ``set of numbers'' with ``intelligible
property of numbers''; see \S \ref{sect2e}).

\noindent 14. Since this statement itself refers to ``all sets'', one
could ask whether a predicativist could accept it without falling into
contradiction, but I suppose he could have a reasonably clear idea of
what is forbidden despite being unable to formally define it. In any case,
the statement is made for the benefit of non-predicativists. A predicativist
should not need to be explicitly forbidden from talking about ``all sets''
since the concept would make no sense to him and it should not even occur
to him to speak this way. A similar comment could apply to the vicious-circle
principle.

\noindent 15. There are several predicatively equivalent versions
of this condition. In intuitionistic logic with arithmetical comprehension,
the numerical omniscience schema, and $\Ac \vee \neg \Ac$ for all atomic
$\Ac$, for any $a \in \omega$ and any ordering on $\omega$ the statement
(1) $(\forall X) {\rm TI}(X,a)$
is equivalent to (2) the assertion that $\{b: b \prec a\}$ has no
proper progressive subsets and also to (3) the assertion that for
all $X$, if there exists $b \prec a$ in $X$ then there is a least such $b$.
Assuming dependent choice for arithmetical formulas, the preceding
are also equivalent to (4) the assertion that every decreasing sequence
in $\{b: b \prec a\}$ is eventually constant and (5) the assertion that
there is no strictly decreasing sequence in $\{b: b \prec a\}$.

\noindent 16. According to reference \cite{Fef2a} it is the $\Delta_0$
collection schema which makes the KP axioms impredicative. Footnote 7
of \cite{Fef2a} refers to \cite{Kre3} for justification of this point,
but the relevant comment in footnote 4 of \cite{Kre3} explicitly locates
impredicativity in the fact that ``the interpretation of the logical
constants, in particular of $\rightarrow$, is classical''. This seems to
imply that if intuitionistic logic were used then the KP axioms would
be predicatively valid, so that weakening the logical axioms from
classical to intuitionistic would render acceptable non-logical axioms
which allow one to access ordinals beyond \Gm{}. Apparently this possibility
was never pursued.

\noindent 17. In an earlier version of this paper I suggested that $T$
could alternatively be interpreted as meaning ``provable in $S$ augmented
by an infinitary $\omega$-rule''. However, this is not helpful because it
is ambiguous about exactly which proof trees would be covered. If we allow
all proof trees that are well-founded in the sense of admitting induction,
then the predicativist should have the same difficulty accepting the
Tarskian implication $T(\ulcorner \Ac\urcorner)\to \Ac$ as he has in
accepting the condition ($*$) discussed in Section 1.4; the Tarskian
hierarchies would then be impredicative in the same way as Feferman's
autonomous systems. To avoid this difficulty we would have to insist
on using proof trees that are well-founded in a strong enough sense to
admit an inductive proof that any theorem proven along such a tree is
actually true. But this requires explicitly using the concept of truth,
which is what the alternative interpretation of $T$ was meant to avoid.
}


\end{document}